\documentclass[a4paper, 10pt]{article}

\usepackage{amsmath, amsthm, amsfonts, latexsym, amscd, amssymb, enumerate, hyperref, xcolor}

\input xypic
\theoremstyle{plain}
\newtheorem{theorem}{Theorem}[section]
\newtheorem{proposition}[theorem]{Proposition} 
\newtheorem{lemma}[theorem]{Lemma}
\newtheorem{corollary}[theorem]{Corollary}
\theoremstyle{remark}
\newtheorem{example}[theorem]{Example}
\newtheorem{remark}[theorem]{Remark}
\theoremstyle{definition}
\newtheorem{definition}[theorem]{Definition}

\title{\textbf{Geometry of Clairaut Riemannian warped product submersions}}

\author{Arkadeepta Roy, Kiran Meena, Hemangi Madhusudan Shah}
\date{}

\begin{document}
	\maketitle
	\begin{abstract}
		\noindent In this paper, we introduce and study the concept of \textit{Clairaut Riemannian warped product submersions} between Riemannian warped product manifolds. By generalizing the notion of Clairaut Riemannian submersions to the setting of Riemannian warped product submersions, we define such submersions via a warping function satisfying a Clairaut relation along geodesics. We establish necessary and sufficient conditions under which a Riemannian warped product submersion satisfies the Clairaut condition, showing that it holds if and only if the girth function defining the Clairaut condition has a horizontal gradient, one component of the fibers is totally geodesic, and the other is totally umbilical with mean curvature vector governed by the warping function. We examine the geometric consequences of this structure, study the harmonicity conditions, and the behavior of the Weyl tensor, etc. Additionally, we illustrate the theory with several non-trivial examples. In the latter part of the paper, we explore a detailed study of the curvature behavior of such submersions. Explicit formulas for the Riemannian, Ricci, and sectional curvature tensors of the source space are derived in terms of the geometry of the target and fiber manifolds, as well as the warping and girth functions. These computations provide geometric insight into how warping and the Clairaut condition affect curvature properties, such as conformal flatness and the non-positivity of certain mixed curvatures. We also analyze the conditions for a trivial warping of the source manifold and for the fibers to be locally symmetric. Furthermore, the Einstein condition has been explored in various scenarios. Finally, we also extend and answer a question posed in \cite{Besse_1987} to the setting of Clairaut warped product submersion. 
	\end{abstract}
	

	\section{Introduction}\label{Intro}
	
	Isometric immersions and Riemannian submersions have been the subject of extensive study and have wide applications, including Yang-Mills theory, Kaluza-Klein theory, supergravity and superstring theories, among others \cite{Falcitelli_2004}. They are also used to construct some Riemannian manifolds with positive or non-negative sectional curvature, as well as Einstein manifolds. These notions become more interesting in the context of product manifolds, as Riemannian warped product manifolds have applications in the construction of Schwarzschild and Robertson-Walker cosmological models and in the identification of new classes of Hamiltonian stationary Lagrangian submanifolds \cite{Bishop_1969, Chen_2008, Neill_1983}. Every Riemannian manifold, hence the warped product manifold, can be embedded in some Euclidean space \cite{Moore_1971, Nash_1956, Erken_2021}. Let $\phi_1: M_1 \to N_1$ and $\phi_2: M_2 \to N_2$ be two smooth maps between Riemannian manifolds, and let $\rho: N_1 \to \mathbb{R}^{+}$ and $f:= \rho \circ \phi_1: M_1 \to \mathbb{R}^{+}$ be two smooth functions. Define a smooth map $\phi := (\phi_1, \phi_2): M_1 \times_{f} M_2 \to N_1 \times_{\rho} N_2$ between warped product manifolds such that $\phi(p_1, p_2) = (\phi_1(p_1), \phi_2(p_2))$. Then we have the following notions. 
	\begin{enumerate}
		\item If $\phi_1$ and $\phi_2$ are isometric immersions, then $\phi$ is also an isometric immersion, namely \textit{warped product isometric immersion}.
		
		\item If $\phi_1$ and $\phi_2$ are Riemannian submersions, then $\phi$ is also a Riemannian submersion, namely \textit{Riemannian warped product submersion}. 
	\end{enumerate} 
	We note that warped product isometric immersions have been significantly explored in the literature (see \cite{Chen_2005, Chen_2017, Chen_2002, Nolker_1996, Tojeiro_2007}). On the other hand, the notion of Riemannian warped product submersion was recently introduced by Erken and Murathan \cite{Erken_2021}. Although Erken et al. studied some properties of such mappings in \cite{Erken_MJOM}, we still need to explore more of the geometry of such submersions in depth. With that inspiration, the present paper investigates various geometric properties and their applications to such submersions.
	
	The Clairaut relation states that for every geodesic $c$ on a surface of revolution $M$, $(e^\psi \circ c) \sin \omega$ is constant, where $e^\psi$ is the distance of a point of $M$ from the axis of rotation and $\omega$ is the angle between the tangent vector of the geodesic and the meridian \cite{Pressley_2010, Aso_1991}. Motivated by the importance of the Clairaut relation and geodesics, Bishop \cite{Bishop_1972} introduced the concept of a Clairaut Riemannian submersion as follows.
	\begin{definition}\label{def_by_Bishop} \cite{Bishop_1972}
		A Riemannian submersion $\phi$ between two Riemannian manifolds $M$ and $N$ is said to be a \textit{Clairaut Riemannian submersion}, if there exists a function $r: M \to \mathbb{R}^{+}$ such that for any geodesic $c$ on $M$, the function $(r \circ c) \sin \omega$ is constant, where at any $t$, $\omega(t)$ is the angle between $\dot{c}(t)$ and the horizontal space at $c(t)$.
	\end{definition}

	Furthermore, Meena and Zawadzki extended this notion to the notion of Clairaut conformal submersion in \cite{Meena_Zawadzki}. However, for a particular dilation, both become the same. In this paper, we generalize the concept of Clairaut Riemannian submersions to the notion of \textit{Clairaut Riemannian warped product submersions}, and explore various geometric properties of such submersions to fill the gap. Clairaut Riemannian warped product maps have been recently explored in \cite{YKS_2025}, which are particular Riemannian warped product maps \cite{MSS}. Very recently, \cite{Hasan_2025} explored whether Riemannian submersions preserve geometric quantities, such as the intermediate Ricci curvature.

	The paper is organized as follows. In Section \ref{prelim}, we recall some basic information that is needed for subsequent sections. Section \ref{sec_crwps} is dedicated to the notion of Clairaut Riemannian warped product submersion with non-trivial examples. It also covers such submersions with some conformally changed metrics. Section \ref{sec_harmonicity} covers harmonic conditions for Clairaut Riemannian warped product submersions. Further, in Section \ref{sec_curvature_relations}, we derive the curvature relations, mainly for the Riemannian, scalar, and Ricci curvature tensors. Finally, Section \ref{sec_implications_curvature_relations} covers various geometric implications and important results, such as symmetry, conformal flatness, trivial warping, Einstein condition, etc.
	
	\section{\textbf{Preliminaries}}\label{prelim}
	
	In this section, we review some key notions and results that will be required for our investigation throughout the paper.
	
	Let $\phi: (M^{m}, g) \to (N^{n}, g')$ be a smooth map between two Riemannian manifolds, and let $\phi_{\ast p}: T_p M \to T_{\phi(p)} N$ be its derivative map at $p$. For each regular value $q \in N$, $\phi^{-1}(q)$ is an $(m - n)$ dimensional submanifold of $M$. The submanifolds $\phi^{-1}(q)$, $q \in N$, are called \textit{fibers of $\phi$}. If we assume $\phi_{*p}$ is surjective for all $p \in M$, then considering $\mathcal{V}_p = \ker {\phi_*}_p$ for any $p \in M$, we obtain an integrable distribution $\mathcal{V}$ corresponding to the foliation of $M$ determined by the fibers of $\phi$ such that $\mathcal{V}_p = T_p\phi^{-1}(q)$, where $\phi(p) = q$. Each $\mathcal{V}_p$ is called the \textit{vertical space} at $p$, $\mathcal{V}$ the \textit{vertical distribution}, and the sections of $\mathcal{V}$ the \textit{vertical vector fields}. At any $p \in M$, we have $T_p M = \mathcal{V}_p \oplus \mathcal{H}_p$; $\mathcal{H}_p$ is called the \textit{horizontal space} at $p$, $\mathcal{H}$ the \textit{horizontal distribution}, and the sections of $\mathcal{H}$ the \textit{horizontal vector fields}. Thus, a vector field on $M$ is called \textit{vertical} if it is always tangent to the fibers. Consequently, a vector field on $M$ is called \textit{horizontal} if it is always orthogonal to the fibers. Moreover, a vector field $X$ on $M$ is called \textit{basic} if $X$ is horizontal and $\phi$-related to a vector field $X'$ on $N$, i.e. $\phi_*(X_p) = X'_{\phi(p)}$ for all $p \in M$. In addition, we denote the projection morphism on the distributions $\ker \phi_*$ and $(\ker \phi_*)^\perp$ by $\mathcal{V}$ and $\mathcal{H}$, respectively. Then we have the following notion of a Riemannian submersion.
	\begin{definition}\label{Riem-Sub}
		Let $\phi: (M^{m}, g) \to (N^{n}, g')$ be a smooth map between two Riemannian manifolds.  
		The map $\phi$ is called a \textit{Riemannian submersion} if it satisfies the following properties:
		\begin{enumerate}[$(i)$]
			\item $\phi$ is onto. 
			\item $(\phi_*)_p$ is a surjective mapping of maximal rank $n$ at any point $p \in M$.
			\item $\phi_\ast$ preserves the lengths of the horizontal vectors.
		\end{enumerate}
	\end{definition}
	O'Neill \cite{Neill_1966} defined the fundamental tensors of a Riemannian submersion $\phi$ defined as above. These tensors are extensively used to study the geometry of Riemannian submersions. They are $(1, 2)$-tensors on $M$, and are given by the following formulae:
	\begin{equation}\label{T-Tensor}
		T(E, F) = T_E F = \mathcal{H} \nabla_{\mathcal{V}E} \mathcal{V}F + \mathcal{V} \nabla_{\mathcal{V}E} \mathcal{H}F,
	\end{equation}
	\begin{equation}\label{A-Tensor}
		A(E, F) = A_E F = \mathcal{V} \nabla_{\mathcal{H}E} \mathcal{H}F + \mathcal{H} \nabla_{\mathcal{H}E} \mathcal{V}F,
	\end{equation}
	for any vector fields $E$ and $F$ on $M$, where $\nabla$ denotes the Levi-Civita connection of $(M, g)$. We also have the following lemmas from \cite{Neill_1966}.
	
	\begin{lemma}\label{Lemma_1}
		For any vertical vectors $U$, $W$ and horizontal vectors $X$, $Y$ on $M$, the tensor fields $T$, $A$ satisfy:
		\begin{enumerate}[$(i)$]
			\item $T_U W = T_W U$,
			\item $A_X Y = -A_Y X = \frac{1}{2} \mathcal{V}[X, Y]$.
		\end{enumerate}
	\end{lemma}
	
	\begin{lemma}
		If $X$, $Y$ are basic vector fields on $M$, $\phi$-related to $X', Y'$ respectively, then:
		\begin{enumerate}[$(i)$]
			\item $g(X, Y) = g'(X', Y') \circ \phi$,
			\item $\mathcal{H}[X, Y]$ is basic, $\phi$-related to $[X', Y']$,
			\item $\mathcal{H}(\nabla_X Y)$ is a basic vector field corresponding to $\nabla'_{X'} Y'$, where $\nabla'$ is the connection on $N$,
			\item for any vertical vector field $U$, $[X, U]$ is vertical.
		\end{enumerate}
		Moreover, if $X$ is basic and $U$ is vertical, then $\mathcal{H}(\nabla_U X) = \mathcal{H}(\nabla_X U) = A_X U$.
	\end{lemma}
	
	\noindent In addition, from (\ref{T-Tensor}) and (\ref{A-Tensor}), we have
	\begin{align}\label{Vertical_and_horizonta_part_of_nabla(V,W)}
		& \nabla_V W = T_V W + \hat{\nabla}_V W, \quad
		\nabla_V X = \mathcal{H} \nabla_V X + T_V X, \, \notag \\
		& \nabla_X V = A_X V + \mathcal{V} \nabla_X V, \quad
		\nabla_X Y = \mathcal{H} \nabla_X Y + A_X Y,
	\end{align}
	for $X, Y \in \Gamma(\mathcal{H})$ and $V, W \in \Gamma(\mathcal{V})$, where $\hat{\nabla}_V W = \mathcal{V} \nabla_V W$. On any fiber $\phi^{-1}(q)$, $q \in N$, $\hat{\nabla}$ coincides with the Levi–Civita connection with respect to the metric induced by $g$ on fiber $\phi^{-1}(q)$. 
	
	Observe that $T$ acts on the fibers as the second fundamental form. Restricted to vertical vector fields, it can be seen that $T = 0$ is equivalent to the condition that the fibers are totally geodesic. A Riemannian submersion is called \textit{Riemannian submersion with totally geodesic fibers} if $T$ vanishes identically.
	
	Let $\{U_1, \ldots, U_{m-n}\}$ be an orthonormal frame of $\mathcal{V}$. Then the horizontal vector field $$ H = \frac{1}{m - n} \sum_{i=1}^{m-n} T_{U_i} U_i$$
	is called the \textit{mean curvature vector field of the fibers}. A Riemannian submersion is called \textit{Riemannian submersion with totally umbilical fibers} if $$T_U W = g(U, W) H$$ for $U, W \in \Gamma(\mathcal{V})$. For any $E \in \Gamma(TM)$, $T_E$ and $A_E$ are skew-symmetric operators on $(\Gamma(TM), g)$ that reverse horizontal and vertical distributions, in the following sense: 
	\begin{align*}
		g(T_E F, G) + g(T_E G, F) = 0 \quad {\rm and} \quad g(A_E F, G) + g(A_E G, F) = 0,
	\end{align*} 
	for any $D, E, G \in \Gamma(TM)$. 
	According to Lemma \ref{Lemma_1}, the horizontal distribution $\mathcal{H}$ is \textit{integrable} if and only if $A = 0$.
	
	We denote the Riemannian curvature tensor of $M$, $N$, the vertical and horizontal distributions by $R$, $R'$, $\Hat{R}$, and $R^*$ respectively. Then we have the following equations that provide curvature relations between them \cite[p. 27-28]{Sahin_2017}: 
	\begin{align}\label{Gauss}
		g(R(U,V)W,F) &= g(\hat{R}(U,V)W,F) - g(T_U F, T_V W) + g(T_V F, T_U W), \\
		g(R(U,V)W,X) &= g((\nabla_U T)_V W, X) - g((\nabla_V T)_U W, X), \\
		g(R(X,Y)Z,H) &= g(R^*(X,Y)Z,H) + 2g(A_Z H, A_X Y) \notag \\ &\quad + g(A_Y H, A_X Z) - g(A_X H, A_Y Z),\\
		g(R(X,Y)Z,V) &= -g((\nabla_Z A)_X Y, V) - g(T_V Z, A_X Y) \notag \\ &\quad - g(A_X Z, T_V Y) + g(A_Y Z, T_V X), \\
		g(R(X,Y)V,W) &= -g((\nabla_V A)_X Y, W) + g((\nabla_W A)_X Y, V) \notag \\ &\quad - g(A_X V, A_Y W) + g(A_X W, A_Y V) \notag \\ &\quad + g(T_V X, T_W Y) - g(T_W X, T_V Y),\\
		g(R(X,V)Y,W) &= -g((\nabla_X T)_V W, Y) - g((\nabla_V A)_X Y, W) \notag \\ &\quad + g(T_V X, T_W Y) - g(A_X V, A_Y W),
	\end{align} 
	where $X,Y,Z,H \in \Gamma(\mathcal{H})$ and $U,V,W,F \in \Gamma(\mathcal{V})$. 
	
	Now, we recall the notion of a warped product manifold.
	\begin{definition}
		Let $(M_1^{m_1}, g_1)$ and $(M_2^{m_2}, g_2)$ be two Riemannian manifolds. Let $f: M_1 \to \mathbb{R}$ be a positive smooth function. Then \textit{warped product} $M_1 \times_f M_2$ of $M_1$ and $M_2$ is the Cartesian product $M_1 \times M_2$ with the metric $g = g_1 + f^2 g_2$.
	\end{definition}
	\noindent More precisely, the Riemannian metric $g$ on $M_1 \times_f M_2$ is defined for vector fields $X, Y$ on $M_1 \times M_2$ by $$ g(X, Y) = g_1(\pi_1^*(X), \pi_1^*(Y)) + f^2(\pi_1(\cdot)) g_2(\pi_2^*(X), \pi_2^*(Y)) $$
	where $\pi_1 : M_1 \times M_2 \to M_1$ and $\pi_2 : M_1 \times M_2 \to M_2$ are projections. We recall that these projections are submersions. In this case, it can be easily seen that the fibers $\{x\} \times M_2 = \pi_1^{-1}(x)$ and the leaves $M_1 \times \{y\} = \pi_2^{-1}(y)$ are Riemannian submanifolds of $M_1 \times_f M_2$. For more details, we refer to \cite{Neill_1983} and \cite{Chen_2017}.
	
	The following lemma describes the Levi-Civita connection on a warped product manifold.
	\begin{lemma}\label{WPConnection}
		Let $M = M_1 \times_f M_2$ be a warped product manifold and $\nabla$, $\nabla^1$, and $\nabla^2$ denote the Levi-Civita connections on $M$, $M_1$, and $M_2$, respectively. If $E_1, F_1$ are vector fields on $M_1$ and $E_2, F_2$ are vector fields on $M_2$, then:
		\begin{enumerate}[$(i)$]
			\item $\nabla_{E_1} F_1$ is the lift of $\nabla^1_{E_1} F_1$,
			\item $\nabla_{E_1} E_2 = \nabla_{E_2} E_1 = \frac{E_1(f)}{f} E_2$,
			\item $\text{nor} (\nabla_{E_2} F_2) = -g(E_2, F_2) (\nabla \ln f)$,
			\item $\text{tan} (\nabla_{E_2} F_2)$ is the lift of $\nabla^2_{E_2} F_2$,
		\end{enumerate}
		where $\nabla f$ denotes the gradient of $f$.
	\end{lemma}
	
	The fundamental tensors associated with Riemannian submersions play a central role in the study of Riemannian warped product submersions, in particular. They give rise to the fundamental equations involving these tensors as follows.
	\begin{lemma}\label{wpstensors}
		$\cite{Erken_2021}$
		Let $\phi= (\phi_1, \phi_2): M=M_1 \times_f M_2 \to N=N_1 \times_\rho N_2$ be a Riemannian warped product submersion between two Riemannian warped product manifolds. If $U_i,V_i \in \Gamma(\mathcal{V}_i)$ and $X_i,Y_i \in \Gamma(\mathcal{H}_i), i = 1, 2$, then we have
		\begin{enumerate}[$(i)$]
			\item $T(U_1, V_1) = T_1(U_1, V_1)$, 
			\item $T(U_1, U_2) = 0$, 
			\item $T(U_2, V_2) = T_2(U_2, V_2) - g_M(U_2, V_2) \, \mathcal{H}(\nabla \ln f)$,
			\item $T(V_1, X_1) = T_1(V_1, X_1), \quad \mathcal{H}(\nabla_{V_1} X_1) = \mathcal{H}_1(\nabla^{1}_{V_1} X_1)$,
			\item $T(V_1, X_2) = 0 = \mathcal{V}(\nabla_{X_2} V_1), \quad A_{X_2} V_1 = (V_1(f)/f) X_2 = \mathcal{H}(\nabla_{V_1} X_2)$,
			\item $T(V_2, X_1) = \mathcal{V}(\nabla_{X_1} V_2) = (X_1(f)/f) V_2, \quad A_{X_1} V_2 = 0 = \mathcal{H}(\nabla_{V_2} X_1)$,
			\item $T(V_2, X_2) = T_2(V_2, X_2) = \mathcal{H}_2(\nabla^{2}_{V_2} X_2)$,
			\item $A(X_1, Y_1) = A_1(X_1, Y_1), \quad \mathcal{H}(\nabla_{X_1} Y_1) = \mathcal{H}_1(\nabla^1_{X_1} Y_1)$,
			\item $\mathcal{H}(\nabla_{X_1} X_2) = (X_1(f)/f) X_2 = \mathcal{H}(\nabla_{X_2} X_1), \quad A_{X_1} X_2 = 0 = A_{X_2} X_1$
			\item $A(X_2, Y_2) = A_2(X_2, Y_2), \quad \quad \mathcal{V}(\nabla \ln f) = 0$, 
			\item $\mathcal{H}(\nabla_{X_2} Y_2) = \mathcal{H}_2(\nabla^2_{X_2} Y_2) - g_M(X_2, Y_2) \, \mathcal{H}(\nabla \ln f)$.    
		\end{enumerate}
	\end{lemma}
	
	\section{Clairaut Riemannian warped product submersions}\label{sec_crwps}
	
	In this section, we define Clairaut Riemannian warped product submersions and discuss some non-trivial examples. First, we \textit{introduce} the notion of Clairaut Riemannian warped product submersion, motivated by Bishop's idea of Clairaut Riemannian submersion [see Definition \ref{def_by_Bishop}].
	\begin{definition}
		A Riemannian warped product submersion $\phi$ between two Riemannian warped product manifolds $M=M_1 {\times}_f M_2$ and $N= N_1 {\times}_\rho N_2$ is said to be a \textit{Clairaut Riemannian warped product submersion}, if there exists a function $r: M \to \mathbb{R}^+$ such that for any geodesic $c$ on $M$, the function $(r \circ c) \sin{\omega}$ is constant, where at any $t$, $\omega(t)$ is the angle between $\Dot{c}(t)$ and the horizontal space at $c(t)$.
	\end{definition}
	The following proposition gives necessary and sufficient conditions for a curve on a warped product manifold to be geodesic. This will be used to prove the main result of this section.
	\begin{proposition}\label{Geodesic-condition}
		Let $\phi : (M=M_1 {\times}_f M_2, g) \to (N=N_1 {\times}_\rho N_2, g')$ be a Riemannian warped product submersion. Let $c: I \to M$ be a regular curve on $M$ such that $U(t) = (U_1(t),U_2(t)) = \mathcal{V}\Dot{c}(t)$ and $X(t) = (X_1(t),X_2(t)) = \mathcal{H}\Dot{c}(t)$, i.e. $X_i \in \Gamma (\mathcal{H}_i)$, $U_i \in \Gamma (\mathcal{V}_i)$, $i = 1,2$. Then $c$ is geodesic on $M$ if and only if
		\begin{align*}
			&\mathcal{H}_1 {\nabla^1}_{X_1}X_1 + \mathcal{H}_1 {\nabla^1}_{U_1}X_1 + A_1(X_1,U_1) + T_1(U_1,U_1) \notag \\& + \mathcal{H}_2 {\nabla^2}_{X_2}X_2 + \mathcal{H}_2 {\nabla^1}_{U_2}X_2 + A_2(X_2,U_2) + T_2(U_2,U_2) \notag \\& + 2 \frac{X_1(f)}{f} X_2 + 2 \frac{U_1(f)}{f} X_2 -(g(U_2,U_2) + g(X_2,X_2)) \mathcal{H}(\nabla \ln{f}) = 0 
		\end{align*}
		and
		\begin{align*}
			&\mathcal{V}_1 {\nabla^1}_{X_1}U_1 + \mathcal{V}_1 {\nabla^1}_{U_1}U_1 + T_1(U_1,X_1) + 2 \frac{X_1(f)}{f} U_2 \notag \\& +
			\mathcal{V}_2 {\nabla^2}_{X_2}U_2 + \mathcal{V}_2 {\nabla^1}_{U_2}U_2 + T_2(U_2,X_2) + 2 \frac{U_1(f)}{f} U_2= 0. 
		\end{align*}
	\end{proposition}
	
	\begin{proof}
		Let $c : I \to M$, $c=(\alpha, \beta)$ be a regular curve on $M=M_1 {\times}_f M_2 $, parametrized by arc length such that $\Dot{c}(t) = U(t) + X(t)$, where $U(t) = \mathcal{V}\Dot{c}(t)$ and $X(t) = \mathcal{H}\Dot{c}(t)$. Let $X = X_1 + X_2$ and $U= U_1 + U_2$, where $X_i \in \Gamma (\mathcal{H}_i)$, $U_i \in \Gamma (\mathcal{V}_i)$, $i = 1,2$. We know that $c$ is a geodesic on $M$ if and only if $\nabla_{\Dot{c}} \Dot{c} = 0$, that is, $$\nabla_{(X_1 + X_2 + U_1 + U_2)} (X_1 + X_2 + U_1 + U_2) = 0.$$ Then, applying Lemma \ref{wpstensors} and separating the horizontal and vertical parts, we obtain the required conditions.
	\end{proof}
	
	Now, we state and prove the main result of this section.
	
	\begin{theorem}(\textbf{The Clairaut condition for Riemannian warped product submersion})\label{clairaut}
		\\Let $\phi= (\phi_1,\phi_2) : (M= M_1 {\times}_f M_2, g) \to (N= N_1 {\times}_\rho N_2, g')$ be a Riemannian warped product submersion with connected fibers. Then $\phi$ is Clairaut with $r=e^{\psi}$, where $\psi : M \to \mathbb{R}$ is a smooth function (called \textit{girth} function), if and only if 
		\begin{enumerate}[$(i)$]
			\item $\nabla \psi$ is horizontal,\label{Result (i) Clairaut}
			\item the fibers of $\phi_1$ are totally umbilical with mean curvature vector field $H_1 = - \nabla \psi|_{M_1} = -\nabla \ln{f}$,\label{Result 2 Clairaut} and
			\item the fibers of $\phi_2$ are totally geodesic.
		\end{enumerate}
	\end{theorem}
	
	\begin{proof}
		Let $c : I \to M$ be a geodesic and for $t\in I$, $\Dot{c} (t) = U (t) + X (t) = U_1 (t) + U_2 (t) + X_1 (t) + X_2 (t)$. Let $\omega(t)$ denote the angle in $[0,\frac{\pi}{2}]$ between $\Dot{c}(t)$ and $X(t)$. Then, similar to the proof of \cite[Theorem 5]{Meena_Zawadzki} and \cite[Theorem 3.2]{Meena_Sahin_Shah}, one can show that $\phi$ is Clairaut with $r=e^\psi$ if and only if we have, along $c$:
		\begin{equation}\label{intermediate-Clairaut}
			g(U,U)g\left(\Dot{c}, {(\nabla\psi)}_{c(t)}\right) + g(T_1(U_1,U_1), X_1) + g(T_2(U_2,U_2), X_2) - \frac{X_1(f)}{f} g(U_2,U_2) = 0,
		\end{equation}
		where $U(t) = (U_1(t),U_2(t)) = \mathcal{V}\Dot{c}(t)$ and $X(t) = (X_1(t),X_2(t)) = \mathcal{H}\Dot{c}(t)$ are vertical and horizontal components of $\Dot{c}(t)$, respectively.
		
		Now we proceed to show that $\nabla \psi$ is horizontal.
		
		\noindent For $\psi : M_1 \times M_2 \to \mathbb{R}$, let $\nabla \psi = W_1 + W_2$, where $W_i \in \Gamma(M_i), i=1,2$. If we consider any geodesic $c : I \to M$ with initial vertical tangent vector, that is, assuming $X = 0$, equivalently, $X_1 = 0, X_2 = 0$, then from (\ref{intermediate-Clairaut}), we have $g(U,U) g(U,\nabla \psi) = 0$, which implies that, $U(\psi) = 0$. Hence, $\psi$ is constant on any fiber, as the fibers are connected. Hence, $\nabla\psi$ is horizontal, that is, $W_1, W_2$ are both horizontal. Thus (\ref{intermediate-Clairaut}) becomes 
		$$ g(U,U)g\left(X, {(\nabla\psi)}_{c(t)}\right) + g(T_1(U_1,U_1), X_1) + g(T_2(U_2,U_2), X_2) - \frac{X_1(f)}{f} g(U_2,U_2) = 0. $$
		Here $g(X, \nabla\psi) = g(X_1 + X_2, W_1 + W_2) = g(X_1, W_1) + g(X_2,W_2)$ and $g(U,U) = g(U_1,U_1) + g(U_2,U_2)$. So we have,
		\begin{align}\label{breaking gradpsi}
			& g(U_1,U_1)g(X_1, W_1) + g(U_1,U_1)g(X_2, W_2) + g(U_2,U_2)g(X_1, W_1) + g(U_2,U_2)g(X_2, W_2) \notag \\ 
			& + g(T_1(U_1,U_1), X_1) + g(T_2(U_2,U_2), X_2) - \frac{X_1(f)}{f} g(U_2,U_2) = 0.
		\end{align}
		
		In the following steps, we establish the remaining statements.
		
		\noindent Let us consider the geodesic with initial tangent vector $ U_1 + X_1 + \lambda (U_2 + X_2) $, for arbitrary $\lambda \neq 0$. Then from (\ref{breaking gradpsi}), we have
		\begin{align}\label{lambda-breaking gradpsi}
			& g(U_1,U_1)g(X_1, W_1) + \lambda g(U_1,U_1)g(X_2, W_2) + {\lambda}^2 g(U_2,U_2)g(X_1, W_1) + {\lambda}^3 g(U_2,U_2)g(X_2, W_2) \notag \\ 
			& + g(T_1(U_1,U_1), X_1) + {\lambda}^3 g(T_2(U_2,U_2), X_2) - {\lambda}^2 \frac{X_1(f)}{f} g(U_2,U_2) = 0.
		\end{align}
		Subtracting (\ref{breaking gradpsi}) and (\ref{lambda-breaking gradpsi}), we obtain
		\begin{align}
			& (1- \lambda) [g(U_1,U_1)g(X_2, W_2) + (1 + \lambda) \left(g(U_2,U_2)g(X_1, W_1) - \frac{X_1(f)}{f} g(U_2,U_2)\right) \notag \\
			& + (1 + \lambda + {\lambda}^2) \left(g(U_2,U_2)g(X_2, W_2) - g(T_2(U_2,U_2), X_2)\right)] = 0, \notag
		\end{align}
		which must hold for all $\lambda \neq 0$, so we must have
		\begin{equation}\label{g(X2,W2)}
			g(U_1,U_1)g(X_2, W_2) = 0,
		\end{equation}
		\begin{equation}\label{g(X1,W1)}
			g(U_2,U_2)g(X_1, W_1) - \frac{X_1(f)}{f} g(U_2,U_2) = 0,
		\end{equation}
		and
		\begin{equation}\label{comparing g(X2,W2)}
			g(U_2,U_2)g(X_2, W_2) - g(T_2(U_2,U_2), X_2) = 0.
		\end{equation}
		Using (\ref{g(X2,W2)}),(\ref{g(X1,W1)}) and (\ref{comparing g(X2,W2)}) in (\ref{breaking gradpsi}), we have
		\begin{equation}\label{comparing g(X1,W1)}
			g(U_1,U_1)g(X_1, W_1) + g(T_1(U_1,U_1), X_1) = 0.
		\end{equation}
		Since we have assumed that $U_1,U_2,X_1,X_2$ are all non-zero, from (\ref{g(X2,W2)}), we affirm $g(X_2, W_2) = 0$. Since $W_2$ is horizontal, we have $W_2 = 0$, thus $\nabla \psi = W_1$. From (\ref{g(X1,W1)}), we conclude $$ g(X_1, W_1) - \frac{X_1(f)}{f} = 0 $$
		\\which gives $$ g(X_1, \nabla \psi - \nabla \ln{f}) = 0. $$
		Since $\nabla \psi$ and $\nabla \ln{f}$ are both horizontal vector fields on $M_1$, we have 
		\begin{equation*}
			\nabla \psi = \nabla \ln{f}.
		\end{equation*}
		From (\ref{comparing g(X2,W2)}), we affirm $$ g(g(U_2,U_2) W_2 - T_2(U_2,U_2), X_2) = 0. $$
		As $W_2 = 0$, we have
		\begin{equation}\label{T2(U2,U2)}
			g(T_2(U_2,U_2), X_2) = 0.
		\end{equation}
		But $T_2(U_2,U_2)$ is a horizontal vector field on $M_2$, we have
		\begin{equation*}
			T_2(U_2,U_2) = 0.
		\end{equation*}
		This shows that $\phi_2$ has totally geodesic fibers, as if the dimension of the fibers is one; this is immediate. And if the fibers are of dimension $\geq 2$, if $V$ and $W$ are orthogonal vertical vectors in $M_2$, then $g(V, W) =0$ and $T(V, W) = T(W, V)$. Then for $U_2 = V + W$ by (\ref{T2(U2,U2)}), we have
		\begin{equation}\label{Polarization 1}
			g(T_2(V,V), X_2) + g(T_2(W,W), X_2) + 2g(T_2(V,W), X_2) = 0.
		\end{equation}
		Similarly, $U_2 = V - W$ in (\ref{T2(U2,U2)}) yields
		\begin{equation}\label{Polarization 2}
			g(T_2(V,V), X_2) + g(T_2(W,W), X_2) - 2g(T_2(V,W), X_2) = 0.
		\end{equation}
		Taking the difference of (\ref{Polarization 1}) and (\ref{Polarization 2}), we have for any two orthogonal vertical vector fields $V$ and $W$ on $M_2$, $T_2(V, W) = 0$. Thus we can write for any vertical $V$, $W$ and horizontal $X_2$ on $M_2$, 
		$$ g(T_2(V,W), X_2) = 0 $$ and hence $ T_2(V,W) = 0 $, which shows that $\phi_2$ has totally geodesic fibers. Again, from (\ref{comparing g(X1,W1)}), we have 
		$$ g(g(U_1,U_1) W_1 + T_1(U_1,U_1), X_1) = 0. $$
		This gives
		$$ T_1(U_1,U_1) = - g(U_1,U_1) W_1 = g(U_1,U_1) (- \nabla \psi). $$
		By the same lines as above, for any vertical vector fields $V$, $W$ on $M_1$, 
		\begin{equation*}
			T_1(V,W) = g(V,W) (- \nabla \psi).
		\end{equation*}
		Consequently, the above equation shows that $\phi_1$ has totally umbilical fibers with mean curvature vector field $$ H_1 = - \nabla \psi.$$
		This proves our theorem.
	\end{proof}

	Now, we have the following immediate corollaries.
	
	\begin{corollary}
		Let $\phi = (\phi_1, \phi_2) : M = M_1 \times_f M_2 \rightarrow N = N_1 \times N_2$ be a Clairaut Riemannian warped product submersion with $r = e^{\psi}$ and horizontal integrable distribution. Let ${L}^{\mathcal{H}_1}, L^{\mathcal{H}_2}$ denote the leaves of the horizontal foliation of $\mathcal{H}_1$ and $\mathcal{H}_2$, and $\mathcal{F}^{\mathcal{V}_1}, \mathcal{F}^{\mathcal{V}_2}$ denote the fibers of $\phi_1$ and $\phi_2$, respectively.
		Then we have the following:
		\begin{enumerate}[$(i)$]
			\item $M_1$ locally splits as a twisted product $ M_1 = {L}^{\mathcal{H}_1} \times_{\psi} \mathcal{F}^{\nu_1}$ and consequently $\widetilde{M}_1$, the universal cover of $M_1$, splits as $\widetilde{M}_1 = \widetilde{{L}^{\mathcal{H}_1}} \times_{\psi} \widetilde{\mathcal{F}^{\nu_1}}$. In addition, $M_2$ locally splits as a product $M_2 = L^{\mathcal{H}_2} \times \mathcal{F}^{\nu_2}$ and consequently $\quad \widetilde{M}_2 = \widetilde{L^{\mathcal{H}_2}} \times \widetilde{\mathcal{F}^{\nu_2}}$.
			Hence, $M$ locally splits as $$ M = \left( {L}^{\mathcal{H}_1} \times_{\psi} \mathcal{F}^{\nu_1} \right) \times_f \left( L^{\mathcal{H}_2} \times \mathcal{F}^{\nu_2} \right).$$
			
			\item If $\operatorname{Hess} (\psi) \equiv 0$, then we have local splitting of $M_1$ as $M_1 = {L}^{\mathcal{H}_1} \times \mathcal{F}^{\nu_1}$.
			Hence, locally $$ M = \left( {L}^{\mathcal{H}_1} \times \mathcal{F}^{\nu_1} \right) \times_f \left( L^{\mathcal{H}_2} \times \mathcal{F}^{\nu_2} \right).$$
		\end{enumerate}
	\end{corollary}
	
	\begin{proof}
		The proof follows from \cite[Proposition 3]{Ponge}. 
	\end{proof}
	
	\begin{corollary}
		Let $\phi= (\phi_1,\phi_2) : (M= M_1 {\times}_f M_2, g) \to (N= N_1 {\times}_\rho N_2, g')$ be a Clairaut Riemannian warped product submersion with connected fibers and $r=e^{\psi}$. Let $\theta$ be a function on $N$ such that for all $x \in M$, we have $\theta(\phi(x)) = 1$, then $\phi= (\phi_1,\phi_2) : (M, g) \to (N, {\theta}^2 g')$ is a Clairaut Riemannian warped product submersion with $r=e^{\psi}$. 
	\end{corollary}
	
	\begin{proof}
		Clearly, for $X,Y \in \mathcal{H}_p$,
		$$(\theta(\phi(p)))^2 g'(\phi_* X, \phi_* Y) = (\theta(\phi(p)))^2 g(X,Y) = g(X,Y)$$ and hence, 
		$$\phi= (\phi_1,\phi_2) : (M, g) \to (N, {\theta}^2 g')$$ is a Riemannian warped product submersion. In fact, it is a Clairaut Riemannian warped product submersion by Theorem \ref{clairaut}.
	\end{proof}
	
	\begin{corollary}
		Let $\phi= (\phi_1,\phi_2) : (M= M_1 {\times}_f M_2, g) \to (N= N_1 {\times}_\rho N_2, g')$ be a Clairaut Riemannian warped product submersion with connected fibers and $r=e^{\psi}$. Let $\theta$ be a positive function on $M$ such that $\mathcal{V}\nabla\theta = 0$, then $\phi= (\phi_1,\phi_2) : (M, {\theta}^2 g) \to (N, g')$ is a Clairaut Riemannian warped product submersion with $r=\theta e^{\psi}$. 
	\end{corollary}

	\begin{corollary}
		Let $\phi= (\phi_1,\phi_2) : (M= M_1 {\times}_f M_2, g) \to (N= N_1 {\times}_\rho N_2, g')$ be a Clairaut Riemannian warped product submersion with connected fibers and $r=e^{\psi}$. Then $\phi= (\phi_1,\phi_2) : (M, e^{-2f}g) \to (N, g')$ is a Clairaut Riemannian warped product submersion with $r=1$. 
	\end{corollary}
	
	\begin{corollary}\label{HLaplacian_Psi}
		Let $\phi= (\phi_1,\phi_2) : (M= M_1 {\times}_f M_2, g) \to (N= N_1 {\times}_\rho N_2, g')$ be a Clairaut Riemannian warped product submersion with $r=e^{\psi}$. Then, we have the Laplacian of $\psi$, $$\Delta \psi = \Delta^{\mathcal{H}_1}\psi.$$
	\end{corollary}
	
	\begin{proof}
		Since $\phi$ is a Clairaut Riemannian warped product submersion, the fibers of $\phi$ are connected, $\nabla \psi$ is horizontal, and $\nabla \psi |_{M_2} = 0$. Now, in some neighborhood of a fixed point $p \in M$, we choose a parallel basis $\{E^1_1, \dots,E^1_{m_1-n_1}, \tilde{E}^1_{m_1-n_1+1},\dots,$ $\tilde{E}^1_{m_1}, E^2_1,\dots,E^2_{m_2-n_2}, \tilde{E}^2_{m_2-n_2+1},\dots, \tilde{E}^2_{m_2}\}$, where $\{E^1_i\}_{i= 1}^{m_1 - n_1} \subset \Gamma(\mathcal{V}_1)$, $\{\tilde{E}^1_j\}_{j= m_1 - n_1 + 1}^{m_1} \subset \Gamma(\mathcal{H}_1)$, $\{E^2_i\}_{i=1}^{m_2 - n_2} \subset \Gamma(\mathcal{V}_2)$ and $\{\tilde{E}^2_j\}_{j=m_2 - n_2 +1}^{m_2} \subset \Gamma(\mathcal{H}_2)$ denote the orthonormal frames of the vertical and horizontal distributions of $\phi_1$ and $\phi_2$, respectively. Then
		\begin{align*}
			\Delta^{\mathcal{V}_1}\psi & = \sum_{i=1}^{m_1-n_1}\operatorname{Hess}^\psi (E^1_i,E^1_i) = \sum_{i=1}^{m_1-n_1} g\Big(\nabla_{E^1_i}\nabla \psi, E^1_i\Big) \\
			& = - \sum_{i=1}^{m_1-n_1} g\Big(\nabla_{E^1_i} E^1_i,\nabla \psi \Big) = 0.
		\end{align*}
		Similarly, we can show that $\Delta^{\mathcal{V}_2}\psi = 0$ and $\Delta^{\mathcal{H}_2}\psi = 0$. This yields $\Delta \psi = \Delta^{\mathcal{H}_1}\psi$.
	\end{proof}
	
	Now, we construct some examples of the Clairaut Riemannian warped product submersions.
	
	\begin{example}
		Let $\phi_2 : (M, g) \to (N, g')$ be a Riemannian submersion with connected and totally geodesic fibers. Let $\theta$ be a function on $N$ and the metric $g_\theta$ on $M$ given by: $$g_\theta(X, Y) = g(X, Y), \quad g_\theta(X, U) = 0, \quad g_\theta(U, V) = e^{-2\theta \circ \phi_2} g(U, V),$$
		for any $X, Y \in \mathcal{X}^\mathcal{H}(M)$, $U, V \in \mathcal{X}^\mathcal{V}(M)$. This makes ${\phi}_1 : (M, g_\theta) \to (N, g')$ a Clairaut Riemannian submersion, where ${\phi}_1 = \phi_2$ on $M$. Similarly to \cite[Example 1.8]{Falcitelli_2004}, it is easy to see that the fibers of ${\phi}_1$ are totally umbilical with the mean curvature vector field $H_1 = \nabla(\theta \circ \phi_2)$. Now consider $\phi= ({\phi}_1, \phi_2) : (M_1{\times}_f M_2, h) \to (N_1 {\times}_\rho N_2, h')$, with $M_1 = (M,g_\theta)$, $M_2 = (M,g)$, $N_1= N_2 = (N,g')$, equipped with warped product metrics $h= g_\theta + f^2 \, g$ and $h' = g' + \rho^2 \, g'$, where $f= e^{\theta \circ \phi_2}$ and $\rho$ is given by $\rho \circ \phi_2 = f$. Then $\phi$ is a Clairaut Riemannian warped product submersion with $r=e^\psi$, where $\nabla \psi|_{M_1} = \nabla(\theta \circ \phi_2)$ and $\nabla \psi|_{M_2} = 0$. This gives a wide class for examples of Clairaut Riemannian warped product submersions, where one can vary $\phi_2$ and get new examples.
	\end{example}
	
	\noindent
	
	\begin{example}
		Consider a map $\phi_1 : \left(\mathbb{R}^4, g_1\right) \to \left(\mathbb{R}^3, g'_1\right)$, where $g_1,g'_1$ both are standard metrics, defined by $$\phi_1(x_1, x_2, x_3, x_4) = \left(\sqrt{x_1^2 + x_2^2}, x_3, x_4\right).$$ Then $\phi_1$ is a Clairaut Riemannian submersion with $r = e^\theta$ for $\theta = \ln \left( \sqrt{x_1^2 + x_2^2} \right)$. The fibers of $\phi_1$ are totally umbilical with mean curvature vector field $$ H_1 = -\left( \frac{x_2}{\sqrt{x_1^2 + x_2^2}} \frac{\partial}{\partial x_1} + \frac{x_1}{\sqrt{x_1^2 + x_2^2}} \frac{\partial}{\partial x_2} \right) = - \nabla \theta.$$ 
		Also consider the projection map $\phi_2 : \left(\mathbb{R}^{n+k}, g_2\right) \to \left(\mathbb{R}^n, g'_2\right)$, where $g_2,g'_2$ both are the standard metrics, defined by
		$$ \phi_2(y_1,y_2,\dots,y_n,y_{n+1},\dots,y_{n+k}) = (y_1,y_2,\dots,y_n). $$ Then $\phi_2$ is a Riemannian submersion with totally geodesic fibers, and hence a Clairaut Riemannian submersion trivially. Now consider $\phi= (\phi_1,\phi_2) : (\mathbb{R}^4{\times}_f \mathbb{R}^{n+k}, g) \to (\mathbb{R}^3 {\times}_\rho \mathbb{R}^n, g')$, with $g = g_1 + f^2 \, g_2$ and $g' = g'_1 + \rho^2 \, g'_2$, where $f$ and $\rho$ are given by $\rho \circ \phi_1 = f$. Then $\phi$ is a Clairaut Riemannian warped product submersion with $r=e^\psi$, where $\nabla \psi|_{M_1} = \nabla \theta$ and $\nabla \psi|_{M_2} = 0$.   
	\end{example}
	
	\noindent
	
	Now we classify all the Clairaut Riemannian warped product submersions from $\mathbb{R}^n$, considered as a warped product manifold. We will use the following notation to describe a bundle map. A fiber bundle $(E, B, \pi, F)$ is indicated as $$ F \to E \xrightarrow{\pi}B$$ where $B$ is the base space, $E$ is the total space, $F$ is the fiber, and $\pi$ is the bundle map of the fiber bundle.
	
	\begin{example}
		Let $M=\mathbb{R}^{m+1}$, $M_1 = \mathbb{R}_+ = [0,\infty)$, $M_2 = \mathbb{S}^m$. Then $M =M_1 \times_f M_2,$ where $f: \mathbb{R}_+ \to \mathbb{R}_+ : r \mapsto r$ is the warping function and the warped product metric is $g=dr^2 + r^2 d\theta^2$. Then any Clairaut Riemannian warped product submersion $\phi= (\phi_1,\phi_2) : (M= M_1 {\times}_f M_2, g) \to (N= N_1 {\times}_\rho N_2, g')$ with $r'=e^{\psi}$ is classified as follows:
		\begin{enumerate}[$(1)$]
			\item $N_1 = [a,b),$ for some $a, b \in \mathbb{R}$ and $\rho : [a,b) \to \mathbb{R}_+$ given by $\rho(y) = \dfrac{y-a}{b-y}.$
			Then, from the relation $\rho \circ \phi_1 = f$, we have $\phi_1: \mathbb{R}_+ \to [a,b)$ defined by $\phi_1(x) = \dfrac{a+bx}{1+x}.$
			
			\item $N_1 = [0,\infty) = \mathbb{R}_+$ and $\rho : \mathbb{R}_+ \to \mathbb{R}_+$ is a diffeomorphism. Then, from the relation $\rho \circ \phi_1 = f$, $\phi_1: \mathbb{R}_+ \to \mathbb{R}_+$ is also a diffeomorphism.
		\end{enumerate}
		In both cases, we have at any $p\in \mathbb{R}_+$, $\mathcal{V}_p = \ker(\phi_{1*p}) = \{0\}$ and $\mathcal{H}_p = (\ker(\phi_{1*p}))^\perp = \{\frac{\partial}{\partial r}\}$, where $\{\frac{\partial}{\partial r}\}$ is the basis of $T_p\mathbb{R}_+ \cong \mathbb{R}$. Also, $$\nabla f = \frac{\partial}{\partial r},\; \mbox{hence,}\; \nabla \psi|_{M_1} = \frac{\nabla f}{f} = \frac{1}{r}\frac{\partial}{\partial r}.$$
		Now, we want $\phi_2: S^m \to N_2$ to be a Riemannian submersion with totally geodesic fibers. Assuming $1 \leq \dim (\text{fiber of } \phi_2) \leq m - 1$, then as a fiber bundle $\phi_2$ is one of the following types, as classified in \cite{REJr}:
		\begin{enumerate}[$(a)$]
			\item $S^1 \to S^{2n+1} \xrightarrow{\phi_2} \mathbb{CP}(n) \quad \text{for } n \geq 2$
			\item $S^3 \to S^{4n+3} \xrightarrow{\phi_2} \mathbb{QP}(n) \quad \text{for } n \geq 2$
			\item $S^1 \to S^3 \xrightarrow{\phi_2} S^2 \left(\frac{1}{2}\right)$
			\item $S^3 \to S^7 \xrightarrow{\phi_2} S^4 \left(\frac{1}{2}\right)$
			\item $S^7 \to S^{15} \xrightarrow{\phi_2} S^8 \left(\frac{1}{2}\right)$
		\end{enumerate}
	\end{example}
	
	\section{Harmonicity}\label{sec_harmonicity}
	
	In this section, we study the harmonicity for Clairaut Riemannian warped product submersions. To obtain harmonic conditions, we need their second fundamental form, the tension field, etc. We recall these concepts here. The \textit{second fundamental form} $\nabla \phi_\ast$ of a map $\phi : (M,g) \to (N,g')$ between two Riemannian manifolds is defined as in \cite{Nore_1966}
	$$ (\nabla \phi_\ast) (X,Y) = \Tilde{\nabla}_X \phi_*Y - \phi_*\nabla_XY $$
	for any local vector fields $X,Y$ on M, where $\nabla$ is the Levi-Civita connection of $M$ and $\tilde{\nabla}$ is the pullback of the connection $\nabla'$ of $N$ to the induced vector bundle $\phi^{-1}(TN)$. Furthermore, the \textit{tension field} $\tau(\phi)$ is defined as the trace of $\nabla \phi_\ast$, that is $$ \tau(\phi) = \sum_{i=1}^m (\nabla \phi_\ast )(E_i,E_i) ,$$ where $\{E_i\}_{1\leq i\leq m}$ is a local orthonormal frame around a point $p \in M$. Moreover, we say that $\phi$ is a {\it harmonic map} if and only if $\tau(\phi)$ vanishes at each point $p \in M$. Using these concepts, we have the following result.
	
	\begin{proposition}\label{harmonic}
		Let $\phi= (\phi_1,\phi_2) : (M= M_1 {\times}_f M_2, g) \to (N= N_1 {\times}_\rho N_2, g')$ be a Clairaut Riemannian warped product submersion with $r=e^{\psi}$. Then $\phi$ is harmonic if and only if either $H_1 = 0 = \nabla \psi$ or $(m_1 - n_1) = (m_2 - n_2) = 0$.
	\end{proposition}
	
	\begin{proof}
		Here $$\tau(\phi) = \operatorname{trace} (\nabla \phi_*) = \sum_{k=1}^{m_1 +m_2} (\nabla \phi_*)(E_k,E_k),$$ where $\{E_1, \dots, E_{m_1+m_2}\}$ is a basis of $(M_1\times_f M_2)$. Hence, $$\tau(\phi) = \tau^{\ker \phi_*}(\phi) + \tau^{(\ker \phi_*)^\perp}(\phi) = \tau^{\mathcal{V}}(\phi) + \tau^{\mathcal{H}}(\phi).$$
		Now,
		\begin{align*}
			\tau^{\mathcal{V}}(\phi) & = \sum_{i=1}^{m_1-n_1}(\nabla \phi_*)(E^1_i,E^1_i) + \sum_{a=1}^{m_2-n_2}(\nabla \phi_*)(E^2_a,E^2_a) \\
			& = \phi_*\left(\sum_{i=1}^{m_1-n_1}T(E^1_i,E^1_i)\right) + \phi_*\left(\sum_{a=1}^{m_2-n_2}T(E^2_a,E^2_a)\right) \\
			& = \phi_*\left(\sum_{i=1}^{m_1-n_1}T_1(E^1_i,E^1_i)\right) + \phi_*\left(\sum_{a=1}^{m_2-n_2}[T_2(E^2_a,E^2_a)-g(E^2_a,E^2_a)\nabla (\ln{f})]\right) \\
			& = \phi_*\left(\sum_{i=1}^{m_1-n_1}g(E^1_i,E^1_i)(-\nabla \psi)\right) + \phi_*\left(\sum_{a=1}^{m_2-n_2}-g(E^2_a,E^2_a)\nabla \psi\right) \\
			& = \phi_{1*}\Big(-(m_1-n_1)\nabla \psi \Big) + \phi_{1*}\Big(-(m_2-n_2)\nabla \psi\Big) \\
			& = \phi_*\Big(-\{(m_1-n_1) + (m_2-n_2)\}\nabla \psi\Big).
		\end{align*}
		Also,
		\begin{align*}
			\tau^{\mathcal{H}}(\phi) & = \sum_{j=m_1-n_1+1}^{m_1}(\nabla \phi_*)(\tilde{E}^1_j,\tilde{E}^1_j) + \sum_{b=m_2-n_2+1}^{m_2}(\nabla \phi_*)(\tilde{E}^2_b,\tilde{E}^2_b) \\
			& = 0.
		\end{align*}
		Consequently, $$\tau(\phi) = \tau^{\mathcal{V}}(\phi) + \tau^{\mathcal{H}}(\phi) = \phi_*(-\{(m_1-n_1) + (m_2-n_2)\}\nabla \psi) + 0. $$ This implies that $\tau(\phi) = 0$ if and only if $\phi_*(-\{(m_1-n_1) + (m_2-n_2)\}\nabla \psi) = 0$, that is, either $(m_1-n_1) + (m_2 -n_2) =0$ or $\nabla \psi = 0$. But $\phi_1$ and $\phi_2$ are Riemannian submersions, so $m_1 - n_1 \geq 0$ and $m_2 - n_2 \geq 0$, and thus $(m_1-n_1) + (m_2 -n_2) =0$ gives $(m_1-n_1) = (m_2 -n_2) =0$.
	\end{proof}
	
	Thus, we conclude that: 
	
	\begin{corollary}
		Let $\phi= (\phi_1,\phi_2) : (M= M_1 {\times}_f M_2, g) \to (N= N_1 {\times}_\rho N_2, g')$ be a Clairaut Riemannian warped product submersion with $r=e^{\psi}$. Then $\phi$ is harmonic if and only if either $M$ is a product or $\dim(M_1) = \dim(N_1)$ and $\dim(M_2) = \dim(N_2)$.   
	\end{corollary}
	
	\section{Curvature relations}\label{sec_curvature_relations}
	In this section, we derive formulae for the Riemann curvatures, Ricci curvatures, and sectional curvatures for the source manifold $M$ of a Clairaut Riemannian warped product submersion $\phi$ with $r=e^{\psi}$. In the sequel, $\phi = (\phi_1, \phi_2) : M = M_1 \times_f M_2 \to N = N_1 \times_\rho N_2$ will denote a Clairaut Riemannian warped product submersion with $r = e^\psi$ and $X_i, Y_i, Z_i, H_i \in \Gamma(\mathcal{H}_i)$, $U_i, V_i, W_i, F_i \in \Gamma(\mathcal{V}_i)$, $1 \leq i \leq 2$.
	In addition, $\{E^1_i | i= 1, \dots, m_1 - n_1\} \subset \Gamma(\mathcal{V}_1)$, $\{\tilde{E}^1_j | j= m_1 - n_1 + 1, \dots, m_1\} \subset \Gamma(\mathcal{H}_1)$, $\{E^2_a | a= 1, \dots, m_2 - n_2\} \subset \Gamma(\mathcal{V}_2)$ and $\{\tilde{E}^2_b | b= m_2 - n_2 +1, \dots, m_2 \} \subset \Gamma(\mathcal{H}_2)$ denote the orthonormal frames of the vertical and horizontal distributions of $\phi_1$ and $\phi_2$, respectively, in some neighborhood of a fixed point $p \in M$.
	
	In what follows, $R$, $R_1$, and $R_2$ denote the Riemannian curvature tensors of $M$, $M_1$, and $M_2$, respectively, and $\widehat{R}$, $\widehat{R_1}$ and $\widehat{R_2}$ denote the Riemannian curvature tensors of fibers of $\phi$, $\phi_1$ and $\phi_2$, respectively. Then, by \cite{Chen_2017}, we have the following lemma.
	
	\begin{lemma}\label{wp-curv}
		Let $M = M_1 \times_f M_2$ be a warped product manifold. Let ${\rm Hess}^f$ denote the Hessian of $f$. If $E_1, F_1,G_1 \in \Gamma(M_1)$ and $E_2, F_2,G_2 \in \Gamma(M_2)$, then:
		\begin{enumerate}[$(i)$]
			\item $R(E_1, F_1) G_1 = R_1(E_1, F_1) G_1,$ 
			\item $R(E_1, F_2) F_1 = \frac{{\rm Hess}^f(E_1, F_1)}{f} F_2,$
			\item $R(E_1, F_1) F_2 = R(F_2, G_2) E_1 = 0,$
			\item $R(E_1, F_2) G_2 = -\frac{g(F_2, G_2)}{f} \nabla_{E_1} (\nabla f),$
			\item $R(E_2, F_2) G_2 = R_2(E_2, F_2) G_2 + \frac{\|\nabla f\|^2}{f^2} \left( g(E_2, G_2) F_2 - g(F_2, G_2) E_2 \right).$
		\end{enumerate}
	\end{lemma}
	
	Using appropriate applications of Lemmas \ref{WPConnection}, \ref{wpstensors}, \ref{wp-curv}, and Clairaut condition of Theorem \ref{clairaut} together with the curvature relations mentioned in Section \ref{prelim}, and performing some straightforward computations, we obtain the following relations.
	\begin{theorem}\label{curv-phi}
		Let $\phi = (\phi_1, \phi_2) : (M = M_1 \times_f M_2, g) \to (N = N_1 \times_\rho N_2, g')$ be a Clairaut Riemannian warped product submersion with $r = e^\psi$. Then for any $X_i, Y_i, Z_i, H_i \in \Gamma(\mathcal{H}_i)$, $U_i, V_i, W_i, F_i \in \Gamma(\mathcal{V}_i)$, $1 \leq i \leq 2$, the following relations hold:
		\begin{enumerate}[$(1)$]
			\item $R(U_1, V_1, W_1, F_1) = \widehat{R_1}(U_1, V_1, W_1, F_1) - \|\nabla \psi\|^2 [g(U_1,F_1)g(V_1,W_1) - g(U_1,W_1)g(V_1,F_1)]$,
			
			\item $R(U_2, V_2, W_2, F_2) = \widehat{R_2}(U_2, V_2, W_2, F_2) - 2 \|\nabla \psi\|^2 [g(U_2,F_2)g(V_2,W_2) - g(U_2,W_2)g(V_2,F_2)] $,
			
			\item $R(U_1, V_1, W_1, X_1) = g(U_1,W_1)g(\nabla_{V_1} \nabla \psi ,X_1) - g(V_1,W_1)g(\nabla_{U_1} \nabla \psi ,X_1)$,
			
			\item $R(U_2, V_2, W_2, X_2) = 0$,
			
			\item $R(U_1, X_1, Y_1, V_1) = - g(U_1,V_1) \operatorname{Hess}^{\psi}(X_1,Y_1) - X_1(\psi) Y_1(\psi)g(U_1,V_1)\\
			+ g(\nabla_{U_1}(A(X_1,Y_1)), V_1) + g(A(X_1,V_1), A(Y_1,U_1))$,
			
			\item $R(U_2, X_2, Y_2, V_2) = g(\nabla_{U_2}(A(X_2,Y_2)), V_2) + g(A(X_2,V_2), A(Y_2,U_2)) - \|\nabla \psi\|^2 \, g(X_2,Y_2)g(U_2,V_2)$,
			
			\item $R(X_1,Y_1,Z_1,H_1) = R^*_1(X_1,Y_1,Z_1,H_1) + 2\, g(A(Z_1,H_1),A(X_1,Y_1)) \\+ g(A(Y_1,H_1),A(X_1,Z_1)) + g(A(X_1,H_1),A(Y_1,Z_1))$,
			
			\item $R(X_2,Y_2,Z_2,H_2) = R^*_2(X_2,Y_2,Z_2,H_2) + 2\, g(A(Z_2,H_2),A(X_2,Y_2)) + g(A(Y_2,H_2),A(X_2,Z_2)) \\+ g(A(X_2,H_2),A(Y_2,Z_2)) + \|\nabla \psi \|^2 \, \Big[ g(X_2,Z_2) g(Y_2,H_2) - g(Y_2,Z_2) g(X_2,H_2) \Big]$,
			
			\item $R(U_1,U_2,V_1,V_2) = g(U_1,V_1) g(U_2,V_2) ||\nabla \psi||^2$,
			
			\item $R(X_1,X_2,Y_1,Y_2) = \frac{1}{f} \operatorname{Hess}^{f}(X_1,Y_1)g(X_2,Y_2)$,
			
			\item $R(U_1,U_2,V_2,X_1) = - \frac{1}{f} \operatorname{Hess}^{f}(U_1,X_1) g(U_2,V_2)$,
			
			\item $R(U_1,U_2,V_2,V_2) = 0 = - \frac{1}{f} g(U_2,V_2)g(\nabla _{U_1} \nabla f, V_2)$,
			
			\item $R(X_1,U_2,V_2,U_1) = - \frac{1}{f} \operatorname{Hess}^{f}(X_1,U_1) g(U_2,V_2)$,
			
			\item $R(X_1,U_2,V_2,Y_1) = - \frac{1}{f} \operatorname{Hess}^{f}(X_1,Y_1) g(U_2,V_2)$,
			
			\item $R(U_1,X_2,Y_2,V_1) = - \, g(U_1,V_1) g(X_2,Y_2) ||\nabla \psi||^2$,
			
			\item $R(U_1,X_2,Y_2,X_1) = - \frac{1}{f} \operatorname{Hess}^{f}(U_1,X_1) g(X_2,Y_2)$,
			
			\item $R(X_1,X_2,Y_2,V_1) = - \frac{1}{f} \operatorname{Hess}^{f}(X_1,V_1) g(X_2,Y_2)$,
			
			\item $R(X_1,X_2,Y_2,Y_1) = - \frac{1}{f} \operatorname{Hess}^{f}(X_1,Y_1) g(X_2,Y_2)$,
			
			\item $R(U_1,U_2,V_1,E_1) = 0$, for any $E_1 \in \Gamma (TM_1)$,
			
			\item $R(U_1,U_2,V_1,X_2) = 0$,
			
			\item $R(X_1,X_2,Y_1,E_1) = 0,$ \mbox{for any} $E_1 \in \Gamma (TM_1)$,
			
			\item $R(X_1,X_2,Y_1,U_2) = 0$,
			
			\item $R(E_2,G_2,E_1,F) = 0,$ \mbox{for any} $E_1 \in \Gamma (TM_1)$, $E_2,G_2 \in \Gamma (TM_2)$ \mbox{and} $F \in \Gamma (TM)$,
			
			\item $R(U_1,U_2,V_2,X_2) = 0$,
			
			\item $R(U_1,U_2,Y_2,V_2) = 0$,
			
			\item $R(U_1,U_2,Y_2,Y_1) = 0$,
			
			\item $R(U_1,U_2,Y_2,Y_2) = 0$,
			
			\item $R(X_1,U_2,V_2,E_2) = 0$, \mbox{for any} $E_2 \in \Gamma (TM_2)$,
			
			\item $R(X_1,U_2,Y_2,E) = 0,$ 
			\mbox{for any} $E \in \Gamma (TM)$.
		\end{enumerate}
		
	\end{theorem}
	
	\noindent
	
	\begin{remark}\label{components-curv-phi}
		From Theorem $\ref{curv-phi}$, we observe that for $U_i,V_i,W_i \in \mathcal{V}_i$, $i=1,2$, we have
		\begin{align*}
			& \mathcal{V} R(U_1, V_1, W_1) = \widehat{R_1}(U_1, V_1, W_1) - \|\nabla \psi\|^2 [g(V_1,W_1)U_1 - g(U_1,W_1)V_1], \\
			& \mathcal{V} R(U_2, V_2, W_2) = \widehat{R_2}(U_2, V_2, W_2) - 2 \|\nabla \psi\|^2 [g(V_2,W_2)U_2 - g(U_2,W_2)V_2],\\
			& \mathcal{H} R(U_1, V_1, W_1) = g(U_1,W_1)\nabla_{V_1} \nabla \psi - g(V_1,W_1)\nabla_{U_1} \nabla \psi, \\
			& \mathcal{H}R(U_2, V_2, W_2) = 0. 
		\end{align*}
	\end{remark}
	
	Let ${\rm sec}$, ${\rm sec}_1$, and ${\rm sec}_2$ denote the sectional curvatures of $M$, $M_1$, and $M_2$, respectively, and $\hat{\rm sec}$, $\hat{\rm sec}_1$ and $\hat{\rm sec}_2$ denote the sectional curvatures of the fibers of $\phi$, $\phi_1$ and $\phi_2$, respectively. Then, by direct computation using Theorem \ref{curv-phi}, the following relations for the sectional curvature can be obtained. 
	\begin{corollary}\label{sec-phi}
		Let $\phi = (\phi_1, \phi_2) : (M = M_1 \times_f M_2, g) \to (N = N_1 \times_\rho N_2, g')$ be a Clairaut Riemannian warped product submersion with $r = e^\psi$. Then for any $X_i, Y_i, Z_i, H_i \in \Gamma(\mathcal{H}_i)$, $U_i, V_i, W_i, F_i \in \Gamma(\mathcal{V}_i)$, $1 \leq i \leq 2$, the following relations hold:
		\begin{enumerate}[$(1)$]
			\item $\operatorname{sec}(U_1, V_1) = \operatorname{sec}_1(U_1, V_1) = \hat{\operatorname{sec}}_1(U_1, V_1) - ||\nabla \psi||^2$,
			
			\item $\operatorname{sec}(U_2, V_2) = \Big[ \operatorname{sec}_2(U_2, V_2) - \|\nabla \psi\|^2\Big]= \Big[ \hat{\operatorname{sec}}_2(U_2, V_2) - 2\, \|\nabla \psi\|^2 \Big]$,
			
			\item $\operatorname{sec}(X_1, Y_1) = \frac{1}{\|X_1 \wedge Y_1\|^2}\Big[R^*_1(X_1,Y_1,Y_1,X_1) - 3 \|A(X_1,Y_1)\|^2\Big]$,
			
			\item $\operatorname{sec}(X_2, Y_2) = \frac{1}{\|X_2 \wedge Y_2\|^2}\Big[ R^*_2(X_2,Y_2,Y_2,X_2) - 3 \|A(X_2,Y_2)\|^2 \Big] - ||\nabla \psi||^2$,
			
			\item $\operatorname{sec}(U_1, X_1) = \frac{- \|U_1\|^2 [\operatorname{Hess}^\psi (X_1,X_1) + (X_1(\psi))^2] + \|A_{X_1} U_1\|^2}{\|U_1\|^2 \|X_1\|^2}$,
			
			\item $\operatorname{sec}(U_2, X_2) = \frac{\|A(X_2, U_2)\|^2  }{\|U_2\|^2 \|X_2\|^2} - \|\nabla \psi\|^2$.
		\end{enumerate}
	\end{corollary} 
	
	Now, we compute the Ricci curvatures using Theorem \ref{curv-phi}. Let $\operatorname{Ric}$, $\operatorname{Ric}_1$, and $\operatorname{Ric}_2$ denote the Ricci curvatures of $M$, $M_1$, and $M_2$, respectively. And $\hat{\operatorname{Ric}}$, $\hat{\operatorname{Ric}_1}$ and $\hat{\operatorname{Ric}_2}$ denote the Ricci curvature of fibers of $\phi$, $\phi_1$ and $\phi_2$, respectively. And $\operatorname{Ric}^{\mathcal{H}_i}(X_i,Y_i) = \operatorname{Ric}_{{\rm range}~\phi_{i*}}(\phi_{i*}X_i,\phi_{i*}Y_i)$, for $i=1,2$. Then, we have the following relations.
	
	\begin{corollary}\label{ric-phi}
		Let $\phi = (\phi_1, \phi_2) : (M = M_1 \times_f M_2, g) \to (N = N_1 \times_\rho N_2, g')$ be a Clairaut Riemannian warped product submersion with $r = e^\psi$. Then for any $X_i, Y_i, Z_i, H_i \in \Gamma(\mathcal{H}_i)$, $U_i, V_i, W_i, F_i \in \Gamma(\mathcal{V}_i)$, $1 \leq i \leq 2$, the following hold:
		\begin{align*}
			(1) \quad \operatorname{Ric}(U_1, V_1) = & \hat{\operatorname{Ric}}_1(U_1, V_1) - (m_1 - n_1 + m_2 ) \|\nabla \psi\|_1^2 g(U_1, V_1) \\ & - g(U_1, V_1) \Delta^{\mathcal{H}_1}(\psi) + \operatorname{trace}^{\mathcal{H}_1} \Big[g(A(.,U_1),A(.,V_1))\Big], \\
			(2) \quad \operatorname{Ric}(U_2, V_2) = & \hat{\operatorname{Ric}}_2(U_2, V_2) + \operatorname{trace}^{\mathcal{H}_2} \Big[g(A(.,U_2),A(.,V_2))\Big] \\ & - \left( \Delta^{\mathcal{H}_1}\psi + (m_1 - n_1 + 2m_2 - n_2 - 1)\|\nabla \psi\|^2 \right) g(U_2, V_2), \\
			(3) \quad \operatorname{Ric}(X_1, Y_1) = & \operatorname{Ric}^{\mathcal{H}_1}(X_1, Y_1) - [m_2 + (m_1 -n_1)] \Big[\operatorname{Hess}^\psi (X_1,Y_1) + X_1(\psi)Y_1(\psi) \Big] \\ & + \operatorname{div}^{\mathcal{V}_1}(A(X_1,Y_1)) - 3 \operatorname{trace}^{\mathcal{H}_1} \Big[g(A(X_1,.),A(Y_1,.))\Big] \\ & + \operatorname{trace}^{\mathcal{V}_1} \Big[g(A(X_1,.),A(Y_1,.))\Big] , \\
			(4) \quad \operatorname{Ric}(X_2, Y_2) = & \operatorname{Ric}^{\mathcal{H}_2}(X_2, Y_2) -\Big[(m_1 - n_1 + m_2)\|\nabla \psi\|^2 + \Delta^{\mathcal{H}_1} \psi \Big] \, g(X_2,Y_2) \\ & + \operatorname{div}^{\mathcal{V}_2}(A(X_2,Y_2)) - 3 \operatorname{trace}^{\mathcal{H}_2} \Big[g(A(X_2,.),A(Y_2,.))\Big] \\ & + \operatorname{trace}^{\mathcal{V}_2} \Big[g(A(X_2,.),A(Y_2,.))\Big].
		\end{align*}
	\end{corollary}

	\begin{proof}
		Recall that $\{E^1_i | i= 1, \dots, m_1 - n_1\} \subset \Gamma(\mathcal{V}_1)$, $\{\tilde{E}^1_j | j= m_1 - n_1 + 1, \dots, m_1\} \subset \Gamma(\mathcal{H}_1)$, $\{E^2_a | a= 1, \dots, m_2 - n_2\} \subset \Gamma(\mathcal{V}_2)$ and $\{\tilde{E}^2_b | b= m_2 - n_2 +1, \dots, m_2 \} \subset \Gamma(\mathcal{H}_2)$ denote the orthonormal frames of the vertical and horizontal distributions of $\phi_1$ and $\phi_2$, respectively, in some neighborhood of a fixed point $p \in M$. Then, 
		\begin{align}\label{Ricci_form}
			\operatorname{Ric}(U_1, V_1) & = \sum_{i=1}^{m_1-n_1} R(E^1_i, U_1, V_1, E^1_i) + \sum_{j=m_1-n_1+1}^{m_1} R(\tilde{E}^1_j, U_1, V_1, \tilde{E}^1_j) \notag \\ 
			& + \sum_{a=1}^{m_2-n_2} R(E^2_a, U_1, V_1, E^2_a) + \sum_{b=m_2-n_2+1}^{m_2} R(\tilde{E}^2_b, U_1, V_1, \tilde{E}^2_b). 
		\end{align}
		Using Theorem \ref{curv-phi}, we compute
		\begin{align}\label{t1}
			\sum_i R(E^1_i, U_1, V_1, E^1_i) & = \sum_i \left[ \hat{R}_1(E^1_i, U_1, V_1, E^1_i) - \|\nabla \psi\|^2 \Big( g(E^1_i, E^1_i) g(U_1, V_1) - g(E^1_i, V_1) g(U_1, E^1_i) \Big) \right] \notag \\
			& = \hat{\operatorname{Ric}}_1(U_1, V_1) - (m_1 - n_1 -1) \|\nabla \psi\|^2 g(U_1, V_1), 
		\end{align}
		\begin{align}\label{t2}
			\sum_j R(\tilde{E}^1_j, U_1, V_1, \tilde{E}^1_j) & = \sum_j R(V_1, \tilde{E}^1_j, \tilde{E}^1_j, U_1) \notag \\ 
			& = \sum_j \Big[ - g(V_1, U_1) \operatorname{Hess}^\psi(\tilde{E}^1_j, \tilde{E}^1_j) - \tilde{E}^1_j(\psi) \tilde{E}^1_j(\psi) g(U_1, V_1) \notag \\ & \quad \quad + g(\nabla_{V_1} (A(\tilde{E}^1_j, \tilde{E}^1_j)), U_1) + g(A_{\tilde{E}^1_j} V_1, A_{\tilde{E}^1_j} U_1) \Big] \notag \\ 
			& = - \Big[ \Delta^{\mathcal{H}_1} \psi + \|\nabla \psi\|^2 \Big] g(U_1, V_1) + \sum_j g(A_{\tilde{E}^1_j} U_1, A_{\tilde{E}^1_j} V_1),
		\end{align}
		\begin{align}\label{t3}
			\sum_a R(E^2_a, U_1, V_1, E^2_a) & = \sum_a R(V_1, E^2_a, E^2_a, U_1) 
			= - \sum_a g(U_1, V_1) g(E^2_a, E^2_a) \|\nabla \psi\|^2 \notag \\ 
			& = -(m_2 - n_2) \|\nabla \psi\|^2 g(U_1, V_1),
		\end{align}
		and
		\begin{align}\label{t4}
			\sum_b R(\tilde{E}^2_b, U_1, V_1, \tilde{E}^2_b) & = \sum_b R(V_1, \tilde{E}^2_b, \tilde{E}^2_b, U_1) 
			= \sum_b - \, g(U_1, V_1) g(\tilde{E}^2_b, \tilde{E}^2_b) \|\nabla \psi\|^2 \notag \\ 
			& = - n_2 \|\nabla \psi\|^2 g(U_1, V_1). 
		\end{align}
		Substituting the values from (\ref{t1}), (\ref{t2}), (\ref{t3}) and (\ref{t4}) in (\ref{Ricci_form}) we obtain
		\begin{align}
			\operatorname{Ric}(U_1, V_1) = & \hat{\operatorname{Ric}}_1(U_1, V_1) + \Big[ (m_1 - n_1 + m_2) \|\nabla \psi\|^2 + \Delta^{\mathcal{H}_1} \psi \Big] g(U_1, V_1) \notag \\ 
			& - \sum_j g(A_{\tilde{E}^1_j} U_1, A_{\tilde{E}^1_j} V_1).
		\end{align}
		This implies $(1)$. Similarly, we can establish $(2), (3)$ and $(4)$. 
	\end{proof}
	
	\section{Geometric implications of curvature relations}\label{sec_implications_curvature_relations}
	This section is dedicated to the applications of the curvature relations obtained in the previous section. We present each case individually in the subsections.
	
	\subsection{Local symmetry}
	In this subsection, we discuss the local symmetry of the fibers of a Clairaut Riemannian warped product submersion $\phi$. We know that a Riemannian manifold is {\it locally symmetric} if and only if $\nabla R \equiv 0$ \cite{Petersen_2016}. Using this fact, we have the following result.
	
	\begin{theorem}\label{LocSym}
		Let $\phi = (\phi_1, \phi_2): (M = M_1 \times_f M_2, g) \to (N = N_1 \times_\rho N_2, g')$ be a Clairaut Riemannian warped product submersion with $r = e^ \psi$ and $\|\nabla \psi\| = 1$. If $M$ is locally symmetric, then the fibers of $\phi_1$ and $\phi_2$ are also locally symmetric. 
	\end{theorem}
	
	\begin{proof}
		Since $M$ is symmetric, we have $(\nabla_E R)(U, V, W) = 0$ for all $E, U, V, W \in \Gamma(\ker\phi_*)$. This implies that
		\[
		\mathcal{V}(\nabla_E R)(U, V, W) = 0 \quad \text{and} \quad \mathcal{H}(\nabla_E R)(U, V, W) = 0.
		\]
		Thus, for any $E_1, U_1, V_1, W_1 \in \Gamma(\ker\phi_{1*})$, we have 
		\begin{align*}
			0 &= \mathcal{V}(\nabla_{E_1} R)(U_1, V_1, W_1) \\
			&= \mathcal{V}\nabla_{E_1}(R(U_1, V_1, W_1)) - \mathcal{V} R(\nabla_{E_1} U_1, V_1, W_1) - \mathcal{V} R(U_1, \nabla_{E_1} V_1, W_1) - \mathcal{V} R(U_1, V_1, \nabla_{E_1} W_1).
		\end{align*}
		Employing Equation (\ref{Vertical_and_horizonta_part_of_nabla(V,W)}), Lemma \ref{WPConnection}, and Theorem \ref{clairaut}, the above equation yields
		\begin{align*}
			0 &= \mathcal{V}\nabla^1_{E_1}(\mathcal{V}R(U_1, V_1, W_1) + \mathcal{H}R(U_1,V_1,W_1)) - \mathcal{V} R(\widehat{\nabla}^1_{E_1} U_1, V_1, W_1) \\
			&\quad + R(\nabla \psi, V_1, W_1)g(E_1, U_1) - \mathcal{V} R(U_1, \widehat{\nabla}^1_{E_1} V_1, W_1) + R(U_1, \nabla \psi, W_1)g(E_1, V_1) \\
			&\quad - \mathcal{V} R(U_1, V_1, \widehat{\nabla}^1_{E_1} W_1) + R(U_1, V_1, \nabla \psi)g(E_1, W_1).
		\end{align*}
		Using Remark \ref{components-curv-phi} and the fact $\|\nabla \psi\|^2 = 1$ into the above equation, we get
		\begin{align*}
			0 = &\widehat{\nabla}^1_{E_1}(\widehat{R}_1(U_1, V_1,W_1)) - \mathcal{V}\nabla^1_{E_1}(g(V_1, W_1)U_1) + \mathcal{V}\nabla^1_{E_1}(g(U_1, W_1)V_1) + T_1(E_1, \mathcal{H}R(U_1, V_1, W_1)) \\
			&\quad - \widehat{R}_1(\widehat{\nabla}^1_{E_1} U_1, V_1, W_1) + g(V_1, W_1)\widehat{\nabla}^1_{E_1} U_1 - g(\widehat{\nabla}^1_{E_1} U_1, W_1)V_1 \\
			&\quad - \widehat{R}_1(U_1, \widehat{\nabla}^1_{E_1} V_1, W_1) + g(\widehat{\nabla}^1_{E_1} V_1, W_1)U_1 - g(U_1, W_1)\widehat{\nabla}^1_{E_1} V_1 \\
			&\quad - \widehat{R}_1(U_1, V_1, \widehat{\nabla}^1_{E_1} W_1) + g(V_1, \widehat{\nabla}^1_{E_1} W_1)U_1 - g(U_1, \widehat{\nabla}^1_{E_1} W_1)V_1. 
		\end{align*}
		Since we have $T_1(E_1, \mathcal{H}R(U_1, V_1, W_1)) = 0$, we conclude that 
		$$(\widehat{\nabla}^1_{E_1} \widehat{R}_1)(U_1, V_1, W_1) = 0.$$
		Hence, the fibers of $\phi_1$ are locally symmetric. In addition, proceeding similarly, using Remark \ref{components-curv-phi}, Lemma \ref{WPConnection}, and Theorem \ref{curv-phi}, we can show that the fibers of $\phi_2$ are also locally symmetric.
	\end{proof}
	
	\begin{corollary}
		In the setting of the above theorem, the fibers of $\phi_1$ and $\phi_2$ are locally symmetric subspaces of $M$. Moreover, if we suppose that the leaves of the horizontal spaces of $\phi_1$ and $\phi_2$ are integrable and complete, then they must also be locally symmetric subspaces of $M$. Consequently, locally,$$ M = (\mathcal{L}^{\mathcal{H}_1} \times_{\psi} \mathcal{F}^{\mathcal{V}_1}) \times_{f} (\mathcal{L}^{\mathcal{H}_2} \times \mathcal{F}^{\mathcal{V}_2}) ,$$ and the universal covering space of $M$ can be written as a warped product $$ \widetilde{M} = (\widetilde{\mathcal{L}}^{\mathcal{H}_1} \times_{\psi} \widetilde{\mathcal{F}}^{\mathcal{V}_1}) \times_{f} (\widetilde{\mathcal{L}}^{\mathcal{H}_2} \times \widetilde{\mathcal{F}}^{\mathcal{V}_2}) .$$ 
	\end{corollary}
	
	\begin{proposition}
		Suppose that the fibers of $\phi_1$ and $\phi_2$ are complete. If $M$ is a space form having sectional curvature $\kappa$, then the fibers of $\phi_1$ and $\phi_2$ are also space forms having sectional curvature $(\kappa+1)$ and $(\kappa+2)$, respectively. In addition, if the horizontal spaces of $\phi_1$ and $\phi_2$ are integrable and complete, they are also space forms having sectional curvatures $\kappa$ and $(\kappa+1)$, respectively.
	\end{proposition}
	
	\begin{corollary}
		Under the hypothesis of the above proposition, additionally assuming that the leaves of the horizontal spaces of $\phi_1$ and $\phi_2$ are integrable and complete, we get to the following classification:
		\begin{enumerate}
			\item[(i)] If $k = 0$, that is, $\widetilde{M}$ is isometric to Euclidean space $\mathbb{R}^{m_1 + m_2}$, then $\widetilde{\mathcal{F}}^{\mathcal{V}_1}$ and $\widetilde{\mathcal{F}}^{\mathcal{V}_2}$ are isometric to $\mathbb{S}^{m_1 - n_1}(1)$ and $\mathbb{S}^{m_2 - n_2}\left( \frac{1}{\sqrt{2}} \right)$ respectively. 
			Also, $\widetilde{\mathcal{L}}^{\mathcal{H}_1}$ and $\widetilde{\mathcal{L}}^{\mathcal{H}_2}$ must be totally geodesic submanifolds of $\mathbb{R}^{m_1 + m_2}$, so they must be isometric to $\mathbb{R}^{n_1}$ and $\mathbb{S}^{n_2}(1)$ respectively. 
			Hence, we get $$ \mathbb{R}^{m_1 + m_2} = \left( \mathbb{R}^{n_1} \times_{\psi} \mathbb{S}^{m_1 - n_1}(1) \right) \times_{f} \left( \mathbb{S}^{n_2}(1) \times \mathbb{S}^{m_2 - n_2}\left(\frac{1}{\sqrt{2}}\right) \right). $$
			\item[(ii)] If $k = 1$, by a similar argument, we get $$ \mathbb{S}^{m_1 + m_2}(1) = \left( \mathbb{S}^{n_1}(1) \times_{\psi} \mathbb{S}^{m_1 - n_1}\left(\frac{1}{\sqrt{2}}\right) \right) \times_{f} \left( \mathbb{S}^{n_2}\left(\frac{1}{\sqrt{2}}\right) \times \mathbb{S}^{m_2 - n_2}\left(\frac{1}{\sqrt{3}}\right) \right).$$ 
			\item[(iii)] If $k = -1$, then we get $$ \mathbb{H}^{m_1 + m_2} = \left( \mathbb{H}^{n_1} \times_{\psi} \mathbb{R}^{m_1 - n_1} \right) \times_{f} \left( \mathbb{R}^{n_2} \times \mathbb{S}^{m_2 - n_2}(1) \right).$$
		\end{enumerate}
		
	\end{corollary}
	
	\subsection{Local conformal flatness}
	
	Now we turn to establish the relationship between the local conformal flatness of $M$ and that of the fibers of $M_1$ and $M_2$. First, we recall that a Riemannian manifold $(M, g)$ is said to be \textit{locally conformally flat} if every point $p \in M$ has a neighborhood that is conformally equivalent to an open subset of Euclidean space. Now, we prove the following result.
	
	\begin{theorem}\label{locConflat}
		Let $\phi = (\phi_1, \phi_2): (M = M_1 \times_f M_2, g) \to (N = N_1 \times_\rho N_2, g')$ be a Clairaut Riemannian warped product submersion with $r = e^ \psi$. Assume that the fibers of $\phi_1$ and $\phi_2$ are of dimension $\geq 4$. If $M$ is locally conformally flat, then the fibers of $\phi_1$ and $\phi_2$ are also locally conformally flat.
	\end{theorem}
	
	\begin{proof}
		To prove our claim, we use an equivalent criterion given in \cite[Theorem $3.2$ (6)]{Kulkarni_1970}, which says that: the local conformal flatness of $(M^n, g), n \geq 4$ is equivalent to: at every point $p \in M$ and for every quadruple of orthogonal vectors $\{e_1, e_2, e_3, e_4\}$, $$ \operatorname{sec}(e_1,e_2) + \operatorname{sec}(e_3,e_4) = \operatorname{sec}(e_1,e_4) + \operatorname{sec}(e_2,e_3). $$
		Let $p=(p_1,p_2) \in M$ and $\{U_1,V_1,W_1,F_1\}\subset \operatorname{ker}(\phi_{1*p_1})$ be four orthogonal vectors in $T_pM$. Since $M$ is locally conformally flat, we have
		\begin{align*}
			\operatorname{sec}(U_1, V_1) + \operatorname{sec}(W_1, F_1) = \operatorname{sec}(U_1, F_1) + \operatorname{sec}(V_1, W_1).
		\end{align*}
		Using Corollary \ref{sec-phi} in the aforementioned equation, we obtain
		$$ \hat{\operatorname{sec}}_1(U_1, V_1) + \hat{\operatorname{sec}}_1(W_1, F_1) = \hat{\operatorname{sec}}_1(U_1, F_1) + \hat{\operatorname{sec}}_1(V_1, W_1).$$
		Note that, as $\phi_1$ is a Riemannian submersion, its fibers are submanifolds
		of $M_1$. Hence, using the preceding criterion, we conclude that the fibers of $\phi_1$ are locally conformally flat. In addition, one can show the local conformal flatness of the fibers of $\phi_2$ along similar lines.
	\end{proof}
	
	We know that a Riemannian manifold $(M^m, g), m \geq 4$ is locally conformally flat if and only if its Weyl tensor is identically zero (for details, see \cite{Lee_book}). Hence, we conclude that:
	
	\begin{corollary}\label{Cor_Weyl}
		Under the hypotheses of Theorem $\ref{locConflat}$, if $M$ is locally conformally flat, then the Weyl tensors of $M$, $\phi_1$, and $\phi_2$ are identically zero. In other words, $M$, $\phi_1$, and $\phi_2$ is Weyl flat.
	\end{corollary}
	
	\subsection{Trivial warping}
	In this subsection, as another application of Theorem \ref{curv-phi}, we discuss some of the cases where the source manifold $M$ of $\phi$ admits trivial warping. In what follows, let $\phi = (\phi_1, \phi_2): (M = M_1 \times_f M_2, g) \to (N = N_1 \times_\rho N_2, g')$ be a Clairaut Riemannian warped product submersion with connected fibers and $r = e^\psi$. Recall that $\dim(M_1) = m_1$, $\dim(M_2) = m_2$, $\dim(N_1) = n_1$, and $\dim(N_2) = n_2$.  
	
	\begin{theorem}\label{sectional>trivial warping}
		If $\phi$ is a Clairaut Riemannian warped product submersion between $(M = M_1 \times_f M_2, g)$ and $(N = N_1 \times_\rho N_2, g')$ with $r = e^ \psi$, then the warping function $f$ is trivial if 
		\begin{enumerate}[$(i)$]
			\item $\psi$ attains maximum (minimum), provided $\operatorname{sec}(U,X) \leq 0 \, \Big(\geq 0 \Big)$ for all $U \in \Gamma(\mathcal{V}_1)$ and $X \in \Gamma(\mathcal{H}_1)$ for $m_1 > n_1$,
			\item $\psi$ attains maximum (minimum), provided $\operatorname{sec}(X_1,X_2) \leq 0 \, \Big(\geq 0 \Big)$ for all $X_i \in \Gamma(\mathcal{H}_i), \, i=1,2$ for $m_1 = n_1$.
		\end{enumerate}
	\end{theorem}
	
	\begin{proof}
		From Theorem \ref{curv-phi}, we have
		\begin{align*}
			R(U_1, X_1, Y_1, V_1) = & - g(U_1, V_1) \, \operatorname{Hess}^\psi(X_1, Y_1) - X_1(\psi) Y_1(\psi) g(U_1, V_1) \\ & + g(\nabla_{U_1}(A_{X_1} Y_1), V_1) + g(A_{X_1} V_1, A_{Y_1} U_1),
		\end{align*}
		where $U_1,V_1 \in \Gamma(\mathcal{V}_1)$ and $X_1,Y_1 \in \Gamma(\mathcal{H}_1)$.
		Consider parallel orthonormal bases $\{E^1_i | i= 1, \dots, m_1 - n_1\} \subset \Gamma(\mathcal{V}_1)$, $\{\tilde{E}^1_j | j= m_1 - n_1 + 1, \dots, m_1\} \subset \Gamma(\mathcal{H}_1)$ of $\Gamma(M_1)$ in a neighborhood of some fixed point $p \in M$. Then the aforementioned equation yields,
		\begin{align*}
			\sum_i \sum_j R(E^1_i, \tilde{E}^1_j, \tilde{E}^1_j, E^1_i) = & - \sum_i g(E^1_i, E^1_i) \sum_j \operatorname{Hess}^\psi(\tilde{E}^1_j, \tilde{E}^1_j) - \sum_i g(E^1_i, E^1_i) \sum_j g(\nabla \psi, \tilde{E}^1_j)^2 \\ & + \sum_i \sum_j g(\nabla_{E^1_i}(A(\tilde{E}^1_j, \tilde{E}^1_j)), E^1_i) + \sum_i \sum_j g(A(\tilde{E}^1_j, E^1_i), A(\tilde{E}^1_j, E^1_i)).
		\end{align*}
		Consequently,
		\begin{align*}
			\sum_{i,j} \operatorname{sec}(E^1_i, \tilde{E}^1_j) = - (m_1 - n_1) \, \Delta^{\mathcal{H}_1} \psi - (m_1 - n_1) \, \|\nabla \psi\|_1^2 + \sum_{i,j} \|A(\tilde{E}^1_j, E^1_i)\|^2. 
		\end{align*}
		We have $\Delta \psi = \Delta^{\mathcal{H}_1}\psi$ [by Corollary \ref{HLaplacian_Psi}] and $A(\tilde{E}^1_j, E^1_i) = \mathcal{H}\nabla_{\tilde{E}^1_j}E^1_i = 0$. Also, recall that $\nabla \psi |_2 = 0 $. Thus, we have 
		\begin{align}\label{sec(e_i,e_j)}
			\sum_{i,j} \operatorname{sec}(E^1_i, \tilde{E}^1_j) = - (m_1 - n_1) \, \Delta \psi - (m_1 - n_1) \, \|\nabla \psi\|^2 .
		\end{align}
		By hypothesis, $\operatorname{sec} (U,X) \leq 0$ for all $U \in \Gamma(\mathcal{V}_1)$ and $X \in \Gamma(\mathcal{H}_1)$,
		then $\operatorname{sec}(E^1_i, \tilde{E}^1_j) \leq 0$ for all $i= 1, \dots , m_1-n_1$ and $j= m_1-n_1, \dots,m_1$. Then from (\ref{sec(e_i,e_j)}), we have 
		\begin{equation}\label{elliptic phi step 1}
			(m_1 - n_1) \,\Big[ \Delta \psi + \|\nabla \psi\|^2 \Big] \geq 0.
		\end{equation}
		Fix a smooth function $\theta: M \to \mathbb{R}$. Then consider the elliptic operator acting on $\operatorname{C^\infty}(M)$ with respect to $\theta$, defined in \cite{Gilbarg_Trudinger} as $ \Delta_ \theta := \Delta - \nabla \theta$. In particular, choosing $\theta = - \psi$, we obtain 
		$$ \Delta _ {- \psi} \psi = \Delta \psi + \|\nabla \psi\|^2.$$
		Then from (\ref{elliptic phi step 1}), we have 
		\begin{equation}\label{elliptic phi step 2}
			(m_1-n_1) \,\Delta _ {- \psi} \psi \geq 0.
		\end{equation}
		Since $\phi_1$ is a Riemannian submersion, we know $m_1-n_1 \geq 0$.
		\begin{enumerate}[$(i)$]
			\item If $(m_1-n_1) > 0$, then from (\ref{elliptic phi step 2}), we have $\Delta _ {- \psi} \psi \geq 0$ (that is, $\psi$ is subharmonic with respect to $\Delta _ {- \psi}$). Then, invoking the strong maximum principle, we have that if $\psi$ attains a maximum, then $\psi$ is constant. 
			
			\item If $m_1 = n_1$, we use Theorem \ref{curv-phi} $(17)$, which states that for any $X_i,Y_i \in \Gamma(\mathcal{H}_i), \, i= 1,2$, we have $$R(X_1,X_2,Y_2,Y_1) = - \frac{1}{f} \operatorname{Hess}^{f}(X_1,Y_1) g(X_2,Y_2).$$ 
			Proceeding similarly to above, we get 
			\begin{align*}
				\sum_{j,b} \operatorname{sec}(\tilde{E}^1_j, \tilde{E}^2_b) =  - n_2 \, \Big[ \Delta \psi + \|\nabla \psi\|^2 \Big].
			\end{align*}
			By hypothesis, $\operatorname{sec}(X_1,X_2) \leq 0 $ for all $X_i \in \Gamma(\mathcal{H}_i), \, i=1,2$, then $\sum_{j,b} \operatorname{sec}(\tilde{E}^1_j, \tilde{E}^2_b) \leq 0 $ for all $j=1, \dots , m_1$ and $b= m_2-n_2 +1 , \dots,m_2$, and consequently we have $\Delta _ {- \psi} \psi \geq 0$, that is, $\psi$ is subharmonic. Then, invoking the strong maximum principle, we see that if $\psi$ attains a maximum, then $\psi$ is constant.
		\end{enumerate}
		Hence, in both cases $M = M_1 \times M_2$.
	\end{proof}
	Thus, we have the following immediate corollary.
	\begin{corollary}
		Within the framework of Theorem $\ref{sectional>trivial warping}$, if $M_1$ is compact, then $M$ has trivial warping and becomes a product manifold. 
	\end{corollary}
	
	\begin{theorem}\label{const sec curv > flatness}
		Let $\phi$ be a Clairaut Riemannian warped product submersion between $(M = M_1 \times_f M_2, g)$ and $(N = N_1 \times_\rho N_2, g')$ with $r = e^ \psi$. Suppose that $M$ has constant sectional curvature $\kappa$ and $M_1$ is compact. Then the warping function $f$ is trivial, and in this case $M$ and consequently, $M_1$ and $M_2$ are flat. 
	\end{theorem}
	
	\begin{proof}
		Without loss of generality, we assume that $\kappa > 0$. If $\kappa < 0$, then the argument follows in a similar way. Since $M_1$ is compact, $\psi$ reaches a maximum and a minimum. Then by Theorem \ref{sectional>trivial warping}, $f$ is trivial and thus $M$ becomes a product. To prove the flatness, we argue as follows: 
		
		\noindent For any $E_i, F_i \in \Gamma(M_i), \, i=1,2$, we have from Lemma \ref{wp-curv},
		\begin{equation}\label{flatness of M}
			R(E_1, F_2, F_1, E_2) = \frac{\operatorname{Hess}^f(E_1, F_1)}{f} g(F_2,E_2).
		\end{equation}
		Since $f$ is constant, $\operatorname{Hess}^f(E_1,F_1) = 0$. Also, $M$ having constant sectional curvature $\kappa$ implies that 
		$$R(E_1, F_2, F_1, E_2) = \kappa \, \Big[ g(E_1,E_2)g(F_2,F_1) - g(E_1,F_1)g(F_2,E_2) \Big].$$
		Then, (\ref{flatness of M}) gives $\kappa \, g(E_1,F_1) \, g(F_2,E_2) = 0$. Thus, $\kappa = 0$ and therefore, $M$ is flat. Moreover, Lemma \ref{wp-curv} ensures that $M_1$ and $M_2$ are flat.\\ In addition, if $\kappa = 0$, then by a similar argument, $f$ is constant. Consequently, $M$, $M_1$ and $M_2$ are flat by virtue of Lemma \ref{wp-curv}.
	\end{proof}
	
	\begin{remark}
		Let $\widetilde{M}$ denote the universal cover of $M$. In the aforementioned setup, we have $\widetilde{M} = \mathbb{R}^{m_1}/{\Gamma_1} \times \mathbb{R}^{m_2} \quad$ or $\quad \widetilde{M} = {\mathbb{R}}^{m_1}/{\Gamma_1} \times {\mathbb{R}}^{m_2}/{\Gamma_2}$, where $\Gamma_i = \text{Iso}({\mathbb{R}}^{m_i}), \, i=1,2$.
	\end{remark}

	\begin{corollary}
		Let $\phi = (\phi_1, \phi_2): M = M_1 \times_f M_2 \to N = N_1 \times_\rho N_2$ be a Clairaut Riemannian warped product submersion with connected fibers. Under the hypothesis of Theorem $\ref{const sec curv > flatness}$, if $M_1$ is compact, then $N$ has non-positive sectional curvature.
	\end{corollary}
	
	\subsection{Einstein condition}
	This subsection is devoted to exploring the geometry of Clairaut warped product submersion when $M$ is Einstein. By \cite{Besse_1987}, we know that a Riemannian manifold $(M, g)$ is Einstein if its Ricci tensor satisfies $\operatorname{Ric} = \lambda g$, where $\lambda$ is a constant. Finally, we conclude the section by investigating a question posed in \cite{Besse_1987} as an extended set-up of Clairaut warped product submersion. We start with the following results.
	
	\begin{theorem}\label{Einstein}
		Let $\phi = (\phi_1, \phi_2): (M = M_1 \times_f M_2, g) \to (N = N_1 \times_\rho N_2, g')$ be a Clairaut Riemannian warped product submersion with connected fibers, integrable horizontal distribution, and $r = e^ \psi$. If $M$ has a constant sectional curvature, then the fibers of $\phi_i$ are Einstein if $ m_i-n_i \geq 3$, for $i=1,2$.
	\end{theorem}
	
	\begin{proof}
		For any $U_1, V_1 \in \Gamma(\mathcal{V}_1)$, from Corollary \ref{ric-phi} we have 
		\begin{align*}
			\hat{\operatorname{Ric}}_1(U_1, V_1) = \operatorname{Ric}(U_1, V_1) + \Big[(m_1 - n_1 + m_2 ) \|\nabla \psi\|^2 + \Delta^{\mathcal{H}_1} \psi \,\Big] g(U_1, V_1) - \sum_j g(A_{\tilde{E}^1_j} U_1, A_{\tilde{E}^1_j} V_1).
		\end{align*}
		Using the fact that $M$ has constant sectional curvature $\operatorname{sec}$, we have $$\operatorname{Ric}(U_1, V_1) = \operatorname{sec} (m_1 + m_2 - 1) g_1(U_1,V_1).$$
		Also, from the assumption that the horizontal distribution is integrable, we have $A = 0$. Thus, from the above equation, we have the following:
		$$ \hat{\operatorname{Ric}}_1(U_1, V_1) = \Big[(m_1 + m_2 - 1) \operatorname{sec} + (m_1 - n_1 + m_2 ) \|\nabla \psi\|^2 + \Delta^{\mathcal{H}_1} \psi \,\Big] g(U_1, V_1).$$
		If dimension of the fibers of $\phi_1$ is $m_1 - n_1 \geq 3$, then invoking Schur's Lemma, we can conclude that the fibers of $\phi_1$ are Einstein and consequently,
		\begin{equation}\label{C1}
			(m_1 + m_2 - 1) \operatorname{sec} + (m_1 - n_1 + m_2 ) \|\nabla \psi\|^2 + \Delta^{\mathcal{H}_1} \psi \, = C_1 \text{ (constant) .}
		\end{equation}
		By a similar argument, the fibers of $\phi_2$ are Einstein if dimension of the fibers of $\phi_2$ is $m_2 - n_2 \geq 3$ and we get, 
		\begin{equation}\label{C2}
			(m_1 + m_2 - 1) \operatorname{sec} + (m_1 - n_1 + 2 m_2 - n_2 - 1 ) \|\nabla \psi\|^2 + \Delta^{\mathcal{H}_1} \psi \, = C_2 \text{ (constant) .}
		\end{equation}
	\end{proof}
	
	\begin{corollary}
		In the same set-up of Theorem $\ref{Einstein}$, we find that $\psi$ is a distance function (that is, $\|\nabla \psi\|^2 = 1$) with constant $\Delta \psi$.
	\end{corollary}
	
	\begin{proof}
		Comparing the equations (\ref{C1}) and (\ref{C2}), we get $\Delta \psi = C$ (constant) and $$ (m_2 - n_2 - 1) \|\nabla \psi\|^2 = C_2 - C_1 \, , $$
		which shows that $\psi$ is a distance function. 
	\end{proof}
	
	\begin{theorem}\label{Einstein implies const scalar curv}
		Let $\phi = (\phi_1, \phi_2): (M = M_1 \times_f M_2, g) \to (N = N_1 \times_\rho N_2, g')$ be a Clairaut Riemannian warped product submersion with $r = e^\psi$ and connected fibers. If $M$ is Einstein, then the fibers of $\phi_1$ and $\phi_2$ have constant scalar curvatures. Moreover, the scalar curvature restricted to the horizontal space of $\phi_2$ is also constant. 
	\end{theorem}
	
	\begin{proof}
		From Corollary \ref{ric-phi}, we have for $U_1,V_1 \in \Gamma(\mathcal{V}_1)$, 
		\begin{align*}
			\hat{\operatorname{Ric}}_1(U_1, V_1) = & \operatorname{Ric}(U_1, V_1) + (m_1 - n_1 + m_2) \|\nabla \psi \|_1^2 g(U_1, V_1) \\ & + g(U_1, V_1) \Delta^{\mathcal{H}_1}(\psi) - \sum_{j=m_1-n_1+1}^{m_1} g(A(\tilde{E}^1_j, U_1), A(\tilde{E}^1_j, V_1)).
		\end{align*}
		If $M$ is Einstein with $\operatorname{Ric} = \lambda g$, that is $\operatorname{Ric}(U_1, V_1) = \lambda g(U_1, V_1)$, then we affirm, 
		\begin{align*}
			\hat{\operatorname{Ric}}_1(U_1, V_1) = & \Big[ \lambda + (m_1 - n_1 + m_2) ||\nabla \psi||_1^2 + \Delta^{\mathcal{H}_1}(\psi) \Big]g(U_1, V_1) \\ & - \sum_{j=m_1-n_1+1}^{m_1} g(A(\tilde{E}^1_j, U_1), A(\tilde{E}^1_j, V_1)).
		\end{align*}
		Taking trace over the basis $\{E^1_i | i=1, \dots ,m_1-n_1\}$ of $\mathcal{V}_1$, we have
		\begin{align*}
			\operatorname{\hat{sec}}_1 = & \Big[ \lambda + (m_1 - n_1 + m_2) ||\nabla \psi||_1^2 + \Delta^{\mathcal{H}_1}(\psi) \Big](m_1-n_1) \\ & - \sum_{i=1}^{m_1-n_1}\sum_{j=m_1-n_1+1}^{m_1} g(A(\tilde{E}^1_j, E^1_i), A(\tilde{E}^1_j, E^1_i)).
		\end{align*}
		Differentiating with respect to $U_1$, we have 
		\begin{align*}
			\nabla_{U_1} \operatorname{\hat{sec}}_1 = & (m_1 - n_1) \nabla_{U_1} \left[ \lambda + (m_1 - n_1 + m_2) \|\nabla \psi\|^2 + \Delta^{\mathcal{H}_1} \psi \right] \\
			& - \sum_{i=1}^{m_1-n_1}\sum_{j=m_1-n_1+1}^{m_1} \nabla_{U_1} \Big[ g(A(\tilde{E}^1_j, E^1_i), A(\tilde{E}^1_j, E^1_i))\Big],
		\end{align*}
		which shows that $\nabla_{U_1} \operatorname{\hat{sec}}_1 = 0$ and hence $\operatorname{\hat{sec}}_1$ is constant. That is, the fibers of $\phi_1$ have constant scalar curvature. Similarly, we can show that the fibers of $\phi_2$ also have constant scalar curvature. 
		
		Now we proceed to prove the last statement.
		
		\noindent From Corollary \ref{ric-phi}, we have for $X_2,Y_2 \in \Gamma(\mathcal{H}_2)$,
		\begin{align*}
			\operatorname{Ric}^{\mathcal{H}_2}(X_2, Y_2) = & \operatorname{Ric}(X_2, Y_2) + \left[ (m_1 - n_1 + m_2)\|\nabla \psi\|^2 + \Delta^{\mathcal{H}_1} \psi \right] g(X_2, Y_2) \\
			& - \operatorname{div}^{\mathcal{V}_2}(A_{X_2} Y_2) + 3 \sum_b g(A_{X_2} \tilde{E}^2_b, A_{Y_2} \tilde{E}^2_b) - \sum_a g(A_{X_2} E^2_a, A_{Y_2} E^2_a).
		\end{align*}
		If $M$ is Einstein with $\operatorname{Ric} = \lambda g$, that is, $\operatorname{Ric}(X_2, Y_2) = \lambda g(X_2, Y_2)$, then
		\begin{align}\label{Ric-H2}
			\operatorname{Ric}^{\mathcal{H}_2}(X_2, Y_2) = & \left[ \lambda + (m_1 - n_1 + m_2)\|\nabla \psi\|^2 + \Delta^{\mathcal{H}_1} \psi \right] g(X_2, Y_2) \notag \\
			& - \operatorname{div}^{\mathcal{V}_2}(A_{X_2} Y_2) + 3 \sum_b g(A_{X_2} \tilde{E}^2_b, A_{Y_2} \tilde{E}^2_b) - \sum_a g(A_{X_2} E^2_a, A_{Y_2} E^2_a). 
		\end{align}
		Taking trace over the basis $\{\tilde{E}^2_b \mid b = m_1 - n_1 + 1, \ldots, m_2\} \subset \mathcal{H}_2$, we have
		\begin{align}\label{Scalar-H2}
			\operatorname{sec}^{\mathcal{H}_2} = & \left[ \lambda + (m_1 - n_1 + m_2)\|\nabla \psi\|^2 + \Delta^{\mathcal{H}_1} \psi \right] n_2 \notag \\ & - \sum_b \operatorname{div}^{\mathcal{V}_2}(A_{\tilde{E}^2_b} \tilde{E}^2_b) + 3 \sum_b g(A_{\tilde{E}^2_b} \tilde{E}^2_b, A_{\tilde{E}^2_b} \tilde{E}^2_b) - \sum_{a, b} g(A_{\tilde{E}^2_b} E^2_a, A_{\tilde{E}^2_b} E^2_a).
		\end{align}
		Differentiating along $X_2 \in \Gamma(\mathcal{H}_2)$, we get 
		\begin{align*}
			\operatorname{div}^{\mathcal{H}_2} \left( \operatorname{Ric}^{\mathcal{H}_2}(X_2) \right) = \nabla_{X_2}^2 \operatorname{sec}^{\mathcal{H}_2} = 0 - \sum_{a, b} \nabla_{X_2}^2 g(A_2(\tilde{E}^2_b, E^2_a), A_2(\tilde{E}^2_b, E^2_a)).
		\end{align*}
		Now, let us take the divergence of (\ref{Ric-H2}). Then $\operatorname{div}^{\mathcal{H}_2} \text{Ric}^{\mathcal{H}_2}(X_2) = 0$, which implies that $\nabla_{X_2}^2 \operatorname{sec}^{\mathcal{H}_2} = 0$ and hence $\operatorname{sec}^{\mathcal{H}_2} = \text{constant}$.   
	\end{proof}
	
	\begin{corollary}
		Let $\phi = (\phi_1, \phi_2): (M = M_1 \times_f M_2, g) \to (N = N_1 \times_\rho N_2, g')$ be a Clairaut Riemannian warped product submersion with $r = e^ \psi$ and connected fibers. If $\psi$ attains the minimum (maximum), then $M = M_1 \times_f M_2$ is a Riemannian product manifold, provided $M$ is Einstein and $m_2 - (m_1-n_1) \leq 0 \Big( \geq 0 \Big)$.
	\end{corollary}
	
	\begin{proof}
		By Theorem \ref{Einstein implies const scalar curv}, $\operatorname{sec}^{\mathcal{H}_2}$ is constant. Then employing (\ref{Scalar-H2}) we obtain, 
		\begin{equation}\label{meu}
			\Big[ \lambda + (m_1 - n_1 + m_2) \|\nabla \psi\|^2 + \Delta^{\mathcal{H}_1} \psi \Big] \, n_2 = \mu (\text{constant}).
		\end{equation}
		Assuming that the leaves of $\phi_1$, that is, $L^{\mathcal{H}_1}$ are integrable and compact, we have
		$$ \left( m_2 + m_1 - n_1 \right) \int_{L^{\mathcal{H}_1}}\|\nabla \psi\|^2 + \int_{L^{\mathcal{H}_1}} \Delta^{\mathcal{H}} \psi = \int_{L^{\mathcal{H}_1}} (\mu - \lambda \, n_2) \, .$$
		Thus, $$ (m_2 + m_1 - n_1) \int_{L^{\mathcal{H}_1}} \|\nabla \psi\|^2 + 0 = (\mu - \lambda \, n_2) \cdot \operatorname{Vol}(L^{\mathcal{H}_1}) $$
		which implies that $$\frac{m_2 + m_1 - n_1}{\operatorname{Vol}(L^{\mathcal{H}_1})} \int_{L^{\mathcal{H}_1}} \|\nabla \psi\|^2 = (\mu - \lambda \, n_2) .$$
		As $L^{\mathcal{H}_1}$ is compact, there exists some point $p=(p_1,p_2) \in M$ with $p_1 \in L^{\mathcal{H}_1}$ such that $\nabla \psi (p) = 0$, which implies $\|\nabla \psi(p)\|^2 = 0 $ and consequently, $\mu = \lambda \, n_2$ . Hence, using Corollary \ref{HLaplacian_Psi}, we have from (\ref{meu}), $$(m_2 + m_1 - n_1 ) \|\nabla \psi\|^2 + \Delta \psi = 0 \,.$$
		Thus, $\Delta \psi \geq 0$ as $m_2 + (m_1-n_1) \geq 0$, that is, $\psi$ is subharmonic. Then, using the fact that $L^{\mathcal{H}_1}$ is compact, we find that $\psi$ attains a maximum, and consequently, we establish that $\psi$ is constant. If $\nabla \psi = 0$, we find that $\psi$ is constant, which implies that $f$ is constant, and hence $M$ admits a trivial warping.
	\end{proof}
	
	\noindent{\bf Question posed in \cite{Besse_1987}:} In order to build new compact Einstein manifolds from the given ones, it was questioned in \cite{Besse_1987}: Does there exist a compact Einstein warped product manifold with a non-constant warping function? Indeed, Proposition 5 of \cite{Kim_Kim} answers this question by constructing a non-trivial compact Einstein warped product space. \\
	In analogy, we can ask: Does there exist a compact Einstein warped product manifold $M = M_1 \times_f M_2$ admitting a Clairaut Riemannian warped product submersion with $r = e^\psi$? We attempt to answer this question in the following theorem. 
	\begin{theorem}\label{Besse_Extenstion}
		Suppose that $({L}^{\mathcal{H}_1},g_1)$ is a manifold and $f$ is a smooth function on ${L}^{\mathcal{H}_1}$ satisfying, for a constant $\lambda \in \mathbb{R}$ and $m_1, m_2, n_1 \in \mathbb{N}$, $$\text{Ric}^{\mathcal{H}_1} = \lambda\, g_1 + \frac{(m_2 + m_1 - n_1)}{f} \, \operatorname{Hess}^f,$$ then $f$ satisfies $$ f \, \Delta f + (m_1 + m_2 - n_1 - 1) \| \nabla f \|^2 + \lambda f^2 = \mu $$ for a constant $\mu \in \mathbb{R}$. Hence, for a compact Einstein space $({L}^{\mathcal{H}_2},g_2)$ of dimension $(m_1 + m_2 - n_1)$ with $$ \operatorname{Ric}^{\mathcal{H}_2} = \mu \, g_2,$$ we can make a compact Einstein warped product space ${L}^{\mathcal{H}_1} \times_f L^{\mathcal{H}_2}$ with $$\text{Ric} = \lambda \, g \, ,$$ where $g = g_1 + f^2 \, g_2$, $g_1, \, g_2$ being the Riemannian metrics on ${L}^{\mathcal{H}_1}$ and $L^{\mathcal{H}_2}$ respectively. Moreover, if we take compact Einstein manifolds $\mathcal{F}^{\mathcal{V}_1}$ of dimension $(m_1 - n_1)$ and $\mathcal{F}^{\mathcal{V}_2}$ of dimension $(m_2 - n_2)$ satisfying $$ \operatorname{Ric}^{\mathcal{V}_1} = \left[ \lambda + \frac{(m_1 - n_1 + m_2 -1)}{f^2} \| \nabla f \|^2 + \frac{1}{f}\Delta^{\mathcal{H}_1} f \right] g $$ and $$ \operatorname{Ric}^{\mathcal{V}_2} = \left[ \lambda + \frac{(m_1 - n_1 + 2m_2 - n_2 - 2)}{f^2} \| \nabla f \|^2 + \frac{1}{f} \Delta^{\mathcal{H}_1} f \right] g, $$ respectively, then $$ M = \left( {L}^{\mathcal{H}_1} \times_f L^{\mathcal{H}_2} \right) \times_f \left( \mathcal{F}^{\mathcal{V}_1} \times \mathcal{F}^{\mathcal{V}_2} \right) \cong \left( {L}^{\mathcal{H}_1} \times_f \mathcal{F}^{\mathcal{V}_1} \right) \times_f \left( {L}^{\mathcal{H}_1} \times \mathcal{F}^{\mathcal{V}_2} \right) $$ is a compact Einstein warped product manifold admitting a Clairaut Riemannian warped product submersion $\phi = (\phi_1, \phi_2)$ with integrable horizontal distribution, whose fibers are $\mathcal{F}^{\mathcal{V}_1}$ and $\mathcal{F}^{\mathcal{V}_2}$ and the leaves of the horizontal spaces are ${L}^{\mathcal{H}_1}$ and $L^{\mathcal{H}_2}$ respectively.
	\end{theorem}
	
	To prove the above theorem, we need to prove the following proposition.
	
	\begin{proposition}
		Let $\phi = (\phi_1, \phi_2): (M = M_1 \times_f M_2, g) \to (N = N_1 \times_\rho N_2, g')$ be a Clairaut Riemannian warped product submersion with $r = e^\psi$ with connected fibers. If $M$ is Einstein, then the following identity holds for any $X_1 \in \Gamma(\mathcal{H}_1)$ : 
		$$ \nabla^1_{X_1} \Delta^{\mathcal{H}_1} \psi + 2 \operatorname{Hess}^{\psi}(X_1, \nabla \psi) = \Big[ \operatorname{div}^{\mathcal{H}_1}(\operatorname{Hess}^\psi + d\psi \otimes d\psi) \Big](X_1).$$
	\end{proposition}
	
	\begin{proof}
		For $X_1,Y_1 \in \Gamma(\mathcal{H}_1)$, from Corollary \ref{ric-phi}, we have
		\begin{align*}
			\operatorname{Ric}^{\mathcal{H}_1}(X_1, Y_1) = & \operatorname{Ric}(X_1, Y_1) + (m_2 + m_1 -n_1) \Big[\operatorname{Hess}^\psi (X_1,Y_1) + X_1(\psi)Y_1(\psi) \Big] \notag \\ & - \operatorname{div}^{\mathcal{V}_1}(A(X_1,Y_1)) - \sum_i g(A_{X_1} E^1_i, A_{Y_1} E^1_i) + 3 \sum_j g(A_{\tilde{E}^1_j} X_1, A_{\tilde{E}^1_j} Y_1).
		\end{align*}
		Since $M$ is Einstein, $\operatorname{Ric}(X_1, Y_1) = \lambda \, g(X_1, Y_1)$.
		Thus, we have 
		\begin{align}\label{H1-Ric}
			\operatorname{Ric}^{\mathcal{H}_1}(X_1, Y_1) = & \lambda \, g(X_1, Y_1) + (m_2 + m_1 -n_1) \Big[\operatorname{Hess}^\psi (X_1,Y_1) + X_1(\psi)Y_1(\psi) \Big] \notag \\ & - \operatorname{div}^{\mathcal{V}_1}(A(X_1,Y_1)) - \sum_i g(A_{X_1} E^1_i, A_{Y_1} E^1_i) + 3 \sum_j g(A_{\tilde{E}^1_j} X_1, A_{\tilde{E}^1_j} Y_1).
		\end{align}
		Taking trace over $ \{ \tilde{E}^1_j \mid j = m_1 - n_1 + 1, \ldots, m_1 \} $, we get
		\begin{align*}
			\operatorname{sec}^{\mathcal{H}_1} = & \lambda n_1 + (m_1 - n_1 + m_2) \,\Big[ \Delta^{\mathcal{H}_1} \psi + \lVert \nabla \psi \rVert^2 \Big] - \sum_j \operatorname{div}^{\mathcal{V}_1}(A_{\tilde{E}^1_j} \tilde{E}^1_j) \\ 
			& - \sum_{i,j} g(A_{\tilde{E}^1_j} E^1_i, A_{\tilde{E}^1_j} E^1_i) + 3 \sum_j g(A_{\tilde{E}^1_j} \tilde{E}^1_j, A_{\tilde{E}^1_j} \tilde{E}^1_j) ,
		\end{align*}
		which reduces to 
		\begin{align}\label{S-H1}
			\operatorname{sec}^{\mathcal{H}_1} = \lambda n_1 + (m_2 + m_1 - n_1)\Big[ \Delta^{\mathcal{H}_1} \psi + \|\nabla \psi\|^2 \Big] - \sum_{i,j} g(A_{\tilde{E}^1_j} E^1_i, A_{\tilde{E}^1_j} E^1_i).
		\end{align}
		Also, we can choose a parallel basis while taking the trace and get
		\begin{align*}
			\operatorname{sec}^{\mathcal{H}_1} = \lambda n_1 + (m_2 + m_1 - n_1)\Big[ \Delta^{\mathcal{H}_1} \psi + \|\nabla \psi\|^2 \Big].
		\end{align*}
		Differentiating (\ref{S-H1}) along $X_1 \in \Gamma(\mathcal{H}_1)$, we get
		\begin{align*}
			\nabla^1_{X_1} \operatorname{sec}^{\mathcal{H}_1} = (m_2 + m_1 - n_1) \left[ \nabla^1_{X_1} \Delta^{\mathcal{H}_1} \psi + \nabla^1_{X_1} \|\nabla \psi\|^2 \right] - \sum_{i,j} 2 \, g(\nabla^1_{X_1}(A_{\tilde{E}^1_j} E^1_i), A_{\tilde{E}^1_j} E^1_i).
		\end{align*}
		Again, using the parallel basis argument, we have 
		\begin{align}\label{div-Ric-H1-1}
			\operatorname{div}^{\mathcal{H}_1} \operatorname{Ric}^{\mathcal{H}_1}(X_1) = \nabla^1_{X_1} \operatorname{sec}^{\mathcal{H}_1} = (m_2 + m_1 - n_1) \Big[ \nabla^1_{X_1} \Delta^{\mathcal{H}_1} \psi + 2 \operatorname{Hess}^{\psi}(X_1, \nabla \psi) \Big ].
		\end{align}
		Also from (\ref{H1-Ric}), we get
		\begin{align*}
			\operatorname{Ric}^{\mathcal{H}_1}(X_1, Y_1) = & \lambda \, g(X_1, Y_1) + (m_2 + m_1 - n_1) \left( d\psi \otimes d\psi + \operatorname{Hess}^\psi \right)(X_1, Y_1) - \operatorname{div}^{\mathcal{V}_1}(A)(X_1, Y_1) \\ & - \sum_i (\nabla_{E^1_i} \otimes \nabla_{E^1_i})(X_1, Y_1) - 3 \sum_j (A_{\tilde{E}^1_j} \circ A_{\tilde{E}^1_j})(X_1, Y_1),
		\end{align*}
		which shows that 
		\begin{align}\label{Op-H1-Ric}
			\operatorname{Ric}^{\mathcal{H}_1} =& \lambda \, g + (m_2 + m_1 - n_1)(\operatorname{Hess}^\psi + d\psi \otimes d\psi) - \operatorname{div}^{\mathcal{V}_1}(A) \notag \\& - \sum_i (\nabla_{E^1_i} \otimes \nabla_{E^1_i}) - 3 \sum_j (A_{\tilde{E}^1_j} \circ A_{\tilde{E}^1_j}).
		\end{align}
		Now we take the divergence of (\ref{Op-H1-Ric}) and get
		\begin{align*}
			\operatorname{div}^{\mathcal{H}_1} \operatorname{Ric}^{\mathcal{H}_1} &= \lambda \, \operatorname{div}^{\mathcal{H}_1} (g) + (m_2 + m_1 - n_1) \Big[ \operatorname{div}^{\mathcal{H}_1}(\operatorname{Hess}^\psi + d\psi \otimes d\psi) \Big] \\ &\quad - \operatorname{div}^{\mathcal{H}_1}(\operatorname{div}^{\mathcal{V}_1}(A)) - \sum_i \operatorname{div}^{\mathcal{H}_1}(\nabla_{E^1_i} \otimes \nabla_{E^1_i}) - 3 \sum_j \operatorname{div}^{\mathcal{H}_1}(A_{\tilde{E}^1_j} \circ A_{\tilde{E}^1_j}).
		\end{align*}
		When acted on $X_1 \in \Gamma(\mathcal{H}_1)$, this reduces to 
		\begin{align}\label{div_Ric_H1-2}
			\operatorname{div}^{\mathcal{H}_1} \operatorname{Ric}^{\mathcal{H}_1}(X_1) = (m_2 + m_1 - n_1) \Big[ \operatorname{div}^{\mathcal{H}_1}(\operatorname{Hess}^\psi + d\psi \otimes d\psi) \Big](X_1).
		\end{align}
		Comparing (\ref{div-Ric-H1-1}) and (\ref{div_Ric_H1-2}), we get 
		\begin{align}\label{Div_Hessian}
			\nabla^1_{X_1} \Delta^{\mathcal{H}_1} \psi + 2 \operatorname{Hess}^{\psi}(X_1, \nabla \psi) = \Big[ \operatorname{div}^{\mathcal{H}_1}(\operatorname{Hess}^\psi + d\psi \otimes d\psi) \Big](X_1).
		\end{align}
	\end{proof}
	
	\noindent
	Now we come to the proof of Theorem \ref{Besse_Extenstion}. We prove the theorem by extending \cite[Proposition 5]{Kim_Kim}, using a similar technique. The required steps follow from Corollary \ref{ric-phi}, Equation (\ref{Div_Hessian}), and the fact that $A \equiv 0$ for a horizontal integrable distribution.
	
	\section*{CRediT authorship contribution statement}
	Conceptualization, methodology, investigation, validation, writing draft, review, editing, and reading have been performed by all the authors of the paper.
	
	\section*{Ethics approval}
	The submitted work is original and not submitted to more than one journal for simultaneous consideration.
	
	\section*{Consent to participate} Not applicable.
	
	\section*{Consent for publication} Not applicable.
	
	\section*{Code availability} Not applicable.
	
	\section*{Funding} Not applicable.
	
	\section*{Declaration of competing interest} The authors have no conflict of interest and no financial interests in this article.
	
	\section*{Acknowledgments}
	Arkadeepta Roy gratefully acknowledges Harish-Chandra Research Institute, Prayagraj, India, for its doctoral research fellowship. In addition, all the authors are thankful to the reviewer for his/her comments towards improvement.
	
	\section*{Dedication} The author, Kiran Meena, dedicates this paper to the welcome of her first child (Madhav Meena).
	
	\section*{Data availability} No data was used for the research described in the article.

		\noindent A. Roy \\
		Department of Mathematics, Harish-Chandra Research Institute, A CI of HBNI,\\ Chhatnag Road, Jhunsi, Prayagraj, Uttar Pradesh-211019, India.\\	
		E-mails: arkadeeptaroy@hri.res.in / royarkadeepta@gmail.com \\
        ORCID: 0009-0008-1004-090X \\
		
		\noindent K. Meena\\
		Department of Mathematics, Indian Institute of Technology Jodhpur,\\ NH 62, Nagaur Road, Karwar, Jodhpur, Rajasthan-342030, India.\\
		E-mail: kirankapishmeena@gmail.com\\
		ORCID: 0000-0002-6959-5853 \\

        \noindent H. M. Shah\\
		Department of Mathematics, Harish-Chandra Research Institute, A CI of HBNI,\\ Chhatnag Road, Jhunsi, Prayagraj, Uttar Pradesh-211019, India.\\	
		E-mail: hemangimshah@hri.res.in \\
        ORCID: 0000-0002-2665-5094\\

	\end{document}